\documentclass[12pt,leqno]{article}
\usepackage{amsfonts}
\usepackage{amsmath, amsthm, amsfonts, amssymb, color}
\usepackage{mathrsfs}
\usepackage{color}
\theoremstyle{definition}

\newcommand{\scr}[1]{\mathscr #1}
\definecolor{wco}{rgb}{0.5,0.2,0.3}

\numberwithin{equation}{section} \theoremstyle{remark}

\newcommand{\ua}{\uparrow}

\title{{\bf   Dimension-free Convergence Rate  in   Sliced Wasserstein  Distance   for  Empirical Measures of Markov Processes }\footnote{The author is supported by    NSFC (12531007).} }
\author{
{\bf Feng-Yu Wang}\\
Center for Applied Mathematics and KL-AAGDM, Tianjin University,   China\\
wangfy@tju.edu.cn
 }
\begin{document}
\allowdisplaybreaks
\def\R{\mathbb R}  \def\ff{\frac} \def\ss{\sqrt} \def\BB{\mathbf
B}
\def\N{\mathbb N} \def\kk{\kappa} \def\m{{\bf m}}
\def\ee{\varepsilon}\def\ddd{D^*}
\def\dd{\delta} \def\DD{\Delta} \def\vv{\varepsilon} \def\rr{\rho}
\def\<{\langle} \def\>{\rangle}
  \def\nn{\nabla} \def\pp{\partial} \def\E{\mathbb E}
\def\d{\text{\rm{d}}} \def\bb{\beta} \def\aa{\alpha} \def\D{\scr D}
  \def\si{\sigma} \def\ess{\text{\rm{ess}}}\def\s{{\bf s}}
\def\beg{\begin} \def\beq{\begin{equation}}  \def\F{\scr F}
\def\Ric{\mathcal Ric} \def\BBess{\text{\rm{Hess}}}
\def\e{\text{\rm{e}}} \def\ua{\underline a} \def\OO{\Omega}  \def\oo{\omega}
 \def\tt{\tilde}\def\[{\lfloor} \def\]{\rfloor}
\def\cut{\text{\rm{cut}}} \def\P{\mathbb P} \def\ifn{I_n(f^{\bigotimes n})}
\def\C{\scr C}      \def\aaa{\mathbf{r}}     \def\r{r}
\def\gap{\text{\rm{gap}}} \def\prr{\pi_{{\bf m},\varrho}}  \def\r{\mathbf r}
\def\Z{\mathbb Z} \def\vrr{\varrho} \def\ll{\lambda}
\def\L{\scr L}\def\Tt{\tt} \def\TT{\tt}\def\II{\mathbb I}
\def\i{{\rm in}}\def\Sect{{\rm Sect}}
\def\M{\mathbb M}\def\Q{\mathbb Q} \def\texto{\text{o}} \def\LL{\Lambda}
\def\Rank{{\rm Rank}} \def\BB{\mathbb B} \def\i{{\rm i}} \def\BBR{\hat{\R}^d}
\def\to{\rightarrow}\def\l{\ell}\def\iint{\int}\def\gg{\gamma}
\def\EE{\scr E} \def\W{\mathbb W}
\def\A{\scr A} \def\Lip{{\rm Lip}}\def\S{\mathbb S}
\def\B{\scr B}\def\Ent{{\rm Ent}} \def\i{{\rm i}}\def\itparallel{{\it\parallel}}
\def\g{{\mathbf g}}\def\Sect{{\mathcal Sec}}\def\T{\mathcal T}\def\BBB{{\bf B}}
\def\f{\mathbf f} \def\g{\mathbf g}\def\BBL{{\bf L}}  \def\BBG{{\mathbb G}}
\def\BBd{{D^E}} \def\BBdP{D^E_\phi} \def\BBdd{{\bf \dd}} \def\BBs{{\bf s}} \def\GA{\scr A}
\def\BBg{{\bf g}}  \def\BBdd{\psi_B} \def\supp{{\rm supp}}\def\div{{\rm div}}
\def\ddiv{{\rm div}}\def\osc{{\bf osc}}\def\1{{\bf 1}}\def\BBD{\mathbb D}\def\GG{\Gamma}
\def\fdd{\stackrel{f.d.d.}{\longrightarrow}}\def\H{\mathbb H}

\maketitle

\begin{abstract}  To derive  dimension-free convergence rates of empirical measures for Markov processes on a Banach space $\BB$, we adopt the sliced Wasserstein distance (SW distance) induced by a probability measure with full support on the unit ball in the dual space $\BB^*$.
This distance refers to the convergence in finite dimensions, and is topologically  equivalent to the  classical Wasserstein distance when   $\BB =\R^d$. With the SW distance  we  derive dimension-free convergence rates  for the empirical measures of ergodic Markov processes on $\BB$, which can be sharp as   illustrated by concrete examples.
  The study provides an efficient  way to simulate infinite-dimensional distributions  using sample trajectories of Markov processes, avoiding $``$the curse of dimensionality" in the  classical Wasserstein distance. The main results apply  to a broad class of infinite-dimensional models, and   are illustrated by  partially dissipative SPDEs in the end of the paper.
 \end{abstract} \noindent


 \vskip 2cm

 \section{Introduction}

The convergence rate of empirical measures in Wasserstein distance have been intensively investigated in recent years for  i.i.d. samples and continuous time Markov processes,
and in both settings the rate  becomes dramatically slow  in high dimensions.   To simulate distributions in high  or infinite dimensions, we study dimension-free convergence rates for empirical measures under the sliced Wasserstein distance.
In the following, we first  recall some known results on the convergence rate in  Wasserstein distance, then recall   the sliced Wasserstein distance arising from  statistics and matching learning.

Let $\scr P(\R^d)$ be the space of all probability measures on $\R^d$. For any $p\in [1,\infty)$, the $p$-Wasserstein space
$$\scr P_p(\R^d):= \bigg\{\mu\in \scr P(\R^d):\ \mu(|\cdot|^p):=\int_{\R^d}|x|^p\nu(\d x)<\infty\bigg\}$$ is a Polish space under the $p$-Wasserstein distance
$$\W_p(\nu_1,\nu_2):=\inf_{\Pi\in \C(\nu_1,\nu_2)} \bigg(\int_{\R^d} |x-y|^p\Pi(\d x,\d y)\bigg)^{\ff 1p},$$
  where $\C(\nu_1,\nu_2)$ is the set of all coupling measures of $\nu_1$ and $\nu_2$.

Let $\{X_n\}_{n\ge 1}$ be i.i.d. random variables on $\R^d$ with distribution $\mu$, and consider the empirical measure
$$\mu_n:=\ff 1 n \sum_{i=1}^n \dd_{X_i},\ \ \ n\ge 1,$$
where $\dd_x$ is the Dirac measure at point $x$. When $\mu $ is the uniform distribution on $[0,1]^d$,  for any $p\in [1,\infty)$ we have
$$\big(\E\big[\W_p(\mu_n,\mu)^p\big]\big)^{\ff 1 p}\sim \beg{cases} n^{-\ff 1 2},\ &\text{if}\ d=1,\\
n^{-\ff 1 2}\ss{\log n},\ &\text{if}\ d=2,\\
n^{-\ff 1 d},\ &\text{if}\ d\ge 3,\end{cases}$$
where for two   sequences $\{A_n\}$ and $\{B_n\}$ of positive numbers, $A_n\sim B_n$ means
$$0<\liminf_{n\to\infty}  \ff{A_n}{B_n}\le   \limsup_{n\to\infty}  \ff{A_n}{B_n} <\infty.$$
This result can be found in \cite{T},  and when  $d=2$ it  is
 known as Ajtai-Koml\'os-Tusn\'ady (AKT) optimal matching theorem  for which the renormalization limit
 $$\lim_{n\to\infty} \ff n{\log n} \E\big[\W_2(\mu_n,\mu)^2\big] =\ff 1 {4\pi}$$
 is derive in \cite{AT}. See  \cite{BL, FG, HGT, Ledoux} and references therein  for results  on   general distributions $\mu$.

On the other hand, let $X_t$ be the elliptic diffusion process on a compact connected Riemannian manifold $M$ without boundary or with   reflecting boundary, and consider the empirical measure
$$\mu_t:=\ff 1 t \int_0^t \dd_{X_s}\d s,\ \ \ t>0.$$
Then $X_t$ has a unique invariant probability measure $\mu$, and according to \cite{WWZ}, for large $t>0$
$$\big(\E\big[\W_p(\mu_t,\mu)^p\big]\big)^{\ff 1 p}\sim \beg{cases} t^{-\ff 1 2},\ &\text{if}\ d\le 3,\\
t^{-\ff 1 2}\ss{\log t},\ &\text{if}\ d=4,\\
t^{-\ff {1}{d-2}},\ &\text{if}\ d\ge 5.\end{cases}$$
 See \cite{TWZ, WZ} for the renormalization formulas for elliptic diffusions on compact manifolds with dimension $d\le 4$, and see \cite{HM, LW,    W21, W22, W23, W25, W26} for the  study of different models of Markov processes, which include   killed diffusion processes, SPDEs, the fractional Brownian motion torus, and general ergodic Markov processes.

The above results indicate that  the convergence rate of empirical measures under Wasserstein distance,   in both i.i.d. and Markovian settings, becomes dramatically slow as the dimension increases.
This feature is known as  the $``$curse of dimensionality" in computational sciences.
To avoid the high-dimensional complexity of the Wasserstein distance,  the $``$sliced Wasserstein" (SW for short) distance and its generalizations have been applied in
  statistics and machining learning.

 The  SW  distance  is  introduced in \cite{Ra} as  the average   Wasserstein distance between  one-dimensional projections of two probability measures on $\R^d$.  More precisely, for any
 $$\theta\in\S :=\big\{\theta\in \R^d:\ |\theta|=1\big\},$$
 consider the projection $\mu^\theta:=\mu\circ \theta^{-1}$ of $\mu\in \scr P(\R^d)$ in the direction $\theta$, where $\theta(x):= \<\theta,x\>$ and
$$(\mu\circ\theta^{-1})(A)= \mu\big(\{x\in\R^d: \<x,\theta\>\in A\}\big) $$ for  $A\in\B(\R),$   the Borel $\si$-field on $\R$.
Let $\ll_0$ be the uniform distribution on $\S$.
Then  for any $p\in [1,\infty)$,  the    $p$-SW distance is defined as
$$ {\W}_{p,\ll_0}(\mu,\nu):=\bigg(\int_{\S} \W_p(\mu^\theta,\nu^\theta)^p\ll_0(\d\theta)\bigg)^{\ff 1 p},\ \ \ \mu,\nu\in \scr P_p(\R^d).$$
This distance     has been extended in   \cite{Bo} by using an absolutely continuous  probability measure $\ll$ on $\S$ to replace $\ll_0$.  The barycenter of the SW distance was   computed in \cite{BON}.    See \cite{K,M} and references therein for modifications of the SW distance and  applications.

 When  $\R^d$ is replaced by    an infinite-dimensional separable Hilbert space $\H$,     the uniform distribution $\ll_0$ on  the unit ball $\S$ does not exist, so that the above mentioned SW distance  is not available.
As a stronger version of the SW distance,   the  max-SW distance has been  introduced in \cite{24,52}, which   also applies to the separable Hilbert space $\H$:
 $$\overline{\W}_p(\mu,\nu):=\sup_{\theta\in \S} \W_p(\mu^\theta,\nu^\theta),\ \ \ \mu,\nu\in \scr P_p(\H),$$
 where $\scr P_p(\H)$ is the set of all probability measures on $\H$ having finite $p$-th moment, and $\S$ is the unit ball in $\H$.
 This distance   has been used in \cite{HRW} to derive dimension-free convergence rates  for the empirical measure
 $\mu_n$ of  i.i.d. random variables on $\H$ with real distribution $\mu$.
  Since  for each $\theta\in\S$, $\mu_n^\theta$ and $\mu^\theta$ are one-dimensional distributions, it is reasonable to ask for
 dimension-free  convergence rate on $\overline{\W}_p(\mu_n,\mu)$. However,  this is  far from trivial.  For instance, letting $\{e_i\}_{1\le i\le d}\subset \S$ be the canonical base of $\R^d$,
 $\{\xi_i^{(n)}:= \W_p(\mu_n^{e_i},\mu^{e_i})\}_{1\le i\le d}$ are i.i.d. random variables such that
   $$\E\big[\overline{\W}_p(\mu_n,\mu)\big]\ge \E\Big[\max_{1\le i\le d}  \xi_i^{(n)}\Big],
 $$  where  the lower bound may essentially depend on $d$     even its oder of $n$ is dimension-free, see \cite{C} for the  blow up rate of  $$  \E\Big[\max_{1\le i\le d} \xi_i^{(n)}\Big]\to\infty\ \text{as} \ d\to\infty.$$
 According to \cite[Theorem 3.5]{HRW}, for $\mu\in \scr P(\R^d)$ with $M_q(\mu):=\mu(|\cdot|^q)<\infty$ for some $q>2p$, there exists a constant $c \in (0,\infty)$ depending on $p,q$ and $M_q(\mu)$ such that
\beq\label{SS}   \E\big[\overline{\W}_p(\mu_n,\mu)^p\big]   \le c  d \big[\log(n+1)\big]^{\ff{p}{q}+\ff 1 {2}}   n^{-\ff 1 {2}},\ \ \ n\ge 1, \end{equation}
where the upper bound blows up as $d\to\infty$.
Nevertheless,    dimension-free  convergence rates
 on $\E\big[\overline{\W}_p(\mu_n,\mu)^p\big]$ have been presented in \cite[Subsection 3.2]{HRW} for $\mu$ having strong enough concentration properties, so that it can be   efficiently approximated by its finite-dimensional marginals for which the estimate \eqref{SS} applies.
The proof of \eqref{SS}  is   based on the uniform ratio bounds for   the projection distributions of $\mu_n $ and $\mu$, see \cite[Theorem 3.6]{HRW} and \cite[Theorem 3]{51}, which,   however, are not available for the empirical measures of continuous time Markov processes.

To derive  dimension-free convergence  rates   for  empirical measures of   Markov processes on a Banach space $(\BB, \|\cdot\|)$,   we consider  the  $p$-SW distance for $p\in [1,\infty):$
\beq\label{PL} \W_{p,\ll}(\mu,\nu):= \bigg(\int_{\S} \W_p(\mu^\theta,\nu^\theta)^p\,\ll(\d\theta)\bigg)^{\ff 1 p},\ \ \ \mu,\nu\in \scr P(\BB),\end{equation}
  where $\scr P(\BB)$ is the set of all probability measures on $\BB$,   $\ll $  is a probability measure  on  the unit ball
  $$\S:=\big\{\theta\in \BB^*:\ \|\theta\|_{\BB^*}=1\big\} $$  in the dual space $(\BB^*,\|\cdot\|_{\BB^*})$,  and
\beq\label{PPO}\mu^\theta:= \mu\circ \theta^{-1}\in \scr P_p(\R),\ \ \ \theta\in \S,  \ \mu\in \scr P(\BB).\end{equation}
Let
$$\scr P_p(\BB):=\big\{\mu\in \scr P(\BB):\ M_p(\mu):=\mu(\|\cdot\|^p)<\infty\big\},\ \ p\in [1,\infty).$$
By  \eqref{PPO} and the triangle inequality of   $\W_p$, for any $\mu,\nu\in \scr P_p(\BB)$,  the function
$$\S\ni\theta\mapsto  \W_p(\mu^\theta,\nu^\theta)\in [0,\infty)$$
  is Lipschitz continuous.  Moreover, for any $\mu,\nu\in \scr P(\BB)$, we may choose $\{\hat\mu_n,\hat\nu_n\}_{n\ge 1}\subset \scr P_p(\BB)$ such that
  $$ \W_p(\mu^\theta,\nu^\theta)=\lim_{n\to\infty}  \W_p(\hat\mu_n^\theta,\hat\nu_n^\theta),\ \ \theta\in\S.$$
  So,   $\W_p(\mu^\theta,\nu^\theta)$ is measurable in $\theta\in \S$ for any $\mu,\nu\in \scr P(\BB)$, hence  $\W_{p,\ll} (\mu,\nu)$  in \eqref{PL} is well-defined.

   To ensure that   $(\scr P_p(\BB),\W_{p,\ll})$ is a metric space, we will take $\ll$ from the space
  $$ \scr P_{0} (\S):= \big\{\text{probability\  measures\  on}\ \S \ \text{with\  full\ support}\big\}. $$
We will show that  $\W_{p,\ll}$ is topologically equivalent to  $\W_p$ when   $\BB=\R^d$,
 and in general it is related to the convergence in finite-dimensional distributions.

  It is easy to derive dimension-free convergence rate of $\W_{p,\ll}(\mu_n,\mu)$ for i.i.d. samples $\{X_n\}_{n\ge 1}$ on $\BB$, since $\{\theta(X_n)\}_{n\ge 1}$ are i.i.d.   one-dimensional random variables for each
  $\theta\in \S$,   and results in \cite{BL} for one-dimensional matching problems can be used   to estimate   $\E[\W_p(\mu_n^\theta,\mu^\theta)^p]$.

   However, the study on  $\W_{p,\ll}(\mu_t,\mu)$ for Markov processes is different, since   the projection  $\theta(X_t)$ for a  Markov process $X_t$ on $\BB$  is not necessarily a one-dimensional Markov process, so that the existing results
   on the empirical measure of Markov processes do not apply to the SW distance.

In Section 2,  we  introduce some   properties of the SW distance  $\W_{p,\ll}$.      In Sections 3,  we use this distance  to derive dimension-free convergence rates for empirical measures of
exponential ergodic   Markov processes on $\BB$. In Section 4, we  make an extension  to more general ergodic Markov processes.
Finally, we apply the main results to partially dissipative  SPDEs.

\section{SW distance  for distributions on Banach space}
 Let  $\W_{p,\ll}$ be defined in \eqref{PL} for $p\in [1,\infty)$ and $\ll\in\scr P_0(\S)$.
Then
$$\W_{p,\ll}(\mu,\nu)\le \W_p(\mu,\nu)<\infty,\ \ \ \mu,\nu\in \scr P_p(\BB),$$
and  the triangle inequality for $\W_{p,\ll}$ follows from that for  $\W_p$ and $\|\cdot\|_{L^p(\ll)}$. Moreover, by Theorem \ref{T1}(1) below,
$\W_{p,\ll}(\mu,\nu)=0$ implies $\mu=\nu$. Hence, $(\scr P_p(\BB),\W_{p,\ll})$ is a metric space.

When $\BB$ is infinite-dimensional, it is natural to compare $\W_{p,\ll}$ with the   convergence in finite-dimensional distributions.
For any $m\in\mathbb N$ and $(\theta_1,\cdots, \theta_m)\in \S^m$, the finite-dimensional distribution of $\mu\in\scr P(\BB)$ induced by $(\theta_1,\cdots,\theta_m)$ is defined as
$$\mu^{\theta_1,\cdots,\theta_m} := \mu\circ (\theta_1,\cdots,\theta_m)^{-1} \in \scr P(\R^m).$$
Let  $\{\mu_n, \mu\}_{n\ge 1}\subset \scr P(\BB)$.  We say that  $\mu_n$ converges to $\mu$ in finite-dimensional distributions,
and denote  $\mu_n\fdd\mu,$  if as $n\to\infty$
$$\mu_n^{\theta_1,\cdots,\theta_m}\to \mu^{\theta_1,\cdots,\theta_m}\ \text{weakly\ for\ any}\ m\in\mathbb N\ \text{and}\ (\theta_1,\cdots,\theta_m)\in\S^m.$$
  We have the following result on $\W_{p,\ll}$.

\beg{thm}\label{T1}  Let $p\in [1,\infty)$ and $ \ll\in \scr P_0(\S).$
\beg{enumerate}  \item[$(1)$]  Let $\{\mu_n,\mu\}_{n\ge 1}\subset \scr P_p(\BB)$ such that
\beq\label{KZ} \lim_{R\to\infty} \sup_{n\ge 1} \mu_n (\|\cdot\|\ge R) =0.\end{equation}
Then $\W_{p,\ll}(\mu_n,\mu)\to 0$ if and only if   $\mu_n\fdd\mu$ and
\beq\label{AB0}   \lim_{N\to \infty}  \sup_{n\ge 1} \int_{\S} \mu_n^\theta \big((|\cdot|^p-N)^+\big) \ll(\d\theta)=   0.  \end{equation}
Moreover, under    $\mu_n\fdd\mu$,  \eqref{AB0} is equivalent to
\beq\label{AB0'}   \lim_{n\to \infty}  \int_{\S} \mu_n^\theta \big(|\cdot|^p\big) \ll(\d\theta)=   \int_\S \mu^\theta(|\cdot|^p)\ll(\d\theta).  \end{equation}
\item[$(2)$]  Let $(\BB,\<\cdot,\cdot\>)$ be a separable Hilbert space. If  $\{\mu_n\}_{n\ge 1}\subset \scr P_p(\BB)$ is a $\W_{p,\ll}$-Cauchy sequence   satisfying
\beq\label{BD} \sup_{n\ge 1} \mu_n(\|\cdot\|^p)<\infty,\end{equation}
then there exists a unique $\mu\in \scr P_p(\BB)$ such that $\W_{p,\ll}(\mu_n,\mu)\to 0$ as $n\to\infty$.
\item[$(3)$] If $\BB=\R^d$ for some $d\in \mathbb N$, then $(\scr P_{p}(\BB), \W_{p,\ll})$ is a Polish space. Moreover, for any $\{\mu_n,\mu\}_{n\ge 1}\subset \scr P_p(\R^d)$,
$\W_{p,\ll}(\mu_n,\mu)\to 0$ if and only if  $\W_p(\mu_n,\mu)\to 0$.
  \end{enumerate}
\end{thm}

\paragraph{Remark 1.1.}
When $\BB$ is   infinite-dimensional,     $\W_{p,\ll}(\mu_n,\mu)\to 0$ is not enough to imply $\mu_n\fdd\mu$, and the metric $\W_{p,\ll}$ is not complete on $\scr P_p(\BB)$. Below we present two examples to confirm these for $(\BB,\<\cdot,\cdot\>)$ being a separable Hilbert space with an orthonormal  basis $\{e_k\}_{k\ge 1}$.

(a) Let $h_k(\theta):=\<\theta, e_k\>$ for $k\ge 1$ and $\theta\in\S$.  Then
$ \sum_{k=1}^\infty h_k^2=1$ holds on $\S$, so that $|h_k|\le 1$ and
$$\lim_{k\to\infty} h_k(\theta)=0,\ \ \ \theta\in \S.$$
By the dominated convergence theorem, for any $\ll\in \scr P_0(\S)$,
$$  \lim_{k\to\infty} \ll(|h_k|^p)=0,\ \ \ \ p\in [1,\infty).$$
We choose a sequence $k_n\uparrow\infty$ as $n\uparrow\infty$ such that
$$  \ll(|h_{k_n}|^p)\le n^{-1-p},\ \ \ n\ge 1.$$
Let
$$x_n:= n e_{k_n},\ \ \ \mu_n:=\dd_{x_{k_n}},\ \ \ \mu:=\dd_0.$$
We obtain
$$\lim_{n\to\infty } \W_{p,\ll} (\mu_n,\mu)^p= \lim_{n\to\infty } n^p  \ll(|h_{k_n}|^p) \le  \lim_{n\to\infty } n^{-1}  =0.$$
On the other hand, let $\i:=\ss{-1}$ and
$$\xi:= \sum_{n=1}^\infty n^{-1} e_{k_n}\in \BB.$$ We have
$$\hat \mu_n(\xi):= \mu_n(\e^{\i \xi})= \e^{\i \<\xi,x_n\>} =\e^\i \ne 1= \hat \mu(\xi),\ \ \ n\ge 1,$$
so that $\mu_n\fdd\mu$ does not hold.

(b) For any $k\ge 1$ and $\ll\in \scr P_0 (\S),$  we have
$|\theta_k|\le 1$,
 $\ll(|\theta_k|^2)>0$   and
$$ \sum_{k=1}^\infty\ll(|\theta_k|^q)\le   \sum_{k=1}^\infty\ll(|\theta_k|^2)  = \int_\S |\theta|^2\ll(\d\theta)= \ll(\S)=1,\ \ q\in [2,\infty).$$
So, we may choose a sequence $k_n\uparrow\infty$ as $n\uparrow\infty$ such that
\beq\label{CW} \ll(|\theta_{k_n}|^q)\le n^{-2q},\ \ \ n\ge 1.\end{equation}
Let
 $$\mu_n:=\dd_{x_n},\ \ \ x_n:= \sum_{l=1}^n  e_{k_l},\ \ \ n\ge 1.$$
 Then $\{\mu_n\}_{n\ge 1}$ is a $\W_{p,\ll}$-Cauchy sequence for $p\in [1,\infty)$, but $\mu_n\fdd \dd_{x}$ with
$$x:=\sum_{n=1}^\infty    e_{k_n}\not\in \BB,$$  so that $\mu_n$ does not convergence in $(\scr P_p(\BB), \W_{p,\ll})$.
Indeed, for  $q:=2\lor p$, by the definition of $\W_{q,\ll}$,  H\"older's inequality and \eqref{CW},    we obtain
\beg{align*} &\W_{q,\ll}(\mu_n,\mu_m)^q =\int_{\S} \W_q(\dd_{\<x_n,\theta\>}, \dd_{\<x_m,\theta\>})^q\ll(\d\theta)
 =\int_\S \bigg|\sum_{ m+1\le l\le n} \<\theta, e_{k_l}\>\bigg|^q\ll(\d\theta) \\
 &\le    \bigg(\int_\S  \sum_{ m+1\le l\le n} \big|\<\theta, e_{k_l}\>\big|^q  l^{2(q-1)}\ll(\d\theta) \bigg) \bigg(\sum_{ m+1\le l\le n} l^{-2}\bigg)^{q-1}\\
 & =   \bigg(  \sum_{ m+1\le l\le n} \ll(|\theta_{k_l}|^q)   l^{2(q-1)}  \bigg) \bigg(\sum_{ m+1\le l\le n} l^{-2}\bigg)^{q-1}
  \le    \bigg(\sum_{l=m+1}^\infty  l^{-2} \bigg)^{q}, \ \ \ \ n>m\ge 1.\end{align*}
 So, the sequence $\{\mu_n\}_{n\ge 1}$ is $\W_{\ll,q}$-Cauchy, and hence also $\W_{\ll,p}$-Cauchy since $q\ge p$.

\beg{proof}[Proof of Theorem $\ref{T1}(1)$]   Recall that  $\mu_n\fdd\mu$ if and only if
\beq\label{F}\lim_{n\to\infty} \hat\mu_n(\xi)= \hat\mu(\xi),\ \ \ \xi\in \BB^*.\end{equation}

(a) Let \eqref{KZ} hold and $\W_{p,\ll}(\mu_n,\mu)\to 0$. We  verify \eqref{F} so that  $\mu_n\fdd\mu$.  If \eqref{F} does not hold, then there exist $\xi\in\BB^*, \vv>0$ and a sequence $\N\ni k_n\uparrow\infty$ as $n\uparrow\infty$ such that
\beq\label{F2} \big|\hat\mu_{k_n}(\xi)-\hat\mu (\xi)\big|\ge \vv,\ \ n\ge 1.\end{equation}
It is clear that $\xi\ne 0$, so that $\theta_\xi:= \ff \xi{\|\xi\|_{*}}\in \S$. By \eqref{KZ}, we find a constant  $R>0$ such that
$$\sup_{n\ge 1} (\mu_n+\mu) \big(\|\cdot\|_\BB\ge R)\le \ff \vv 8.$$
Let
$$D_\vv:=\bigg\{\theta\in \S:\ \|\theta-\theta_\xi\|_*<\ff \vv{4 R\|\xi\|_* }\bigg\}.$$
Noting that
$$ |\e^{\i \|\xi\|_* \theta(x)} -\e^{\i  \xi(x)}|\le \min\big\{2, \ \|\xi\|_*\| \theta - \theta_\xi \|_* \|x\|_\BB\big\},\ \ \ x\in\BB,$$
we obtain
\beq\label{KZ2}\beg{split}& \sup_{n\ge 1} (\mu+\mu_{n })(|\e^{\i \|\xi\|_* \theta} -\e^{\i  \xi}|) \le  \sup_{n\ge 1} (\mu+\mu_{n })(1_{\{\|\cdot\|<R\}} |\e^{\i \|\xi\|_* \theta} -\e^{\i  \xi}| + 2 1_{\{\|\cdot\|_\BB\ge R\}})  \\
&  \le   R    \|\xi\|_*\| \theta - \theta_\xi \|_*   + 2 \sup_{n\ge 1} (\mu+\mu_{n })(\|\cdot\|_\BB\ge R)
 \le \ff \vv 4 + \ff\vv 4= \ff { \vv}2,\ \ \ \theta\in D_\vv.\end{split}\end{equation}
Since $\ll$ has full support on $\S$, we    have $\ll(D_\vv)>0$.
Moreover,     $\|\theta\|_*=\|\theta_\xi\|_*=1$ implies
\beg{align*} &\big|\e^{\i \|\xi\|_*  \theta (x)}- \e^{\i \|\xi\|_*  \theta (y)}\big|\le \|\xi\|_*|\theta(x)-\theta(y)|.\end{align*}
By Kantorovich's dual formula,  we obtain
$$\W_p(\mu_n^\theta,\mu^\theta)\ge \W_1(\mu_n^\theta,\mu^\theta)\ge \ff 1 {\|\xi\|_*} \big|\mu_n(\e^{\i \|\xi\|_* \theta})- \mu(\e^{\i \|\xi\|_* \theta})\big|.$$
 Combining these with \eqref{F2} and  \eqref{KZ2},  we derive
\beq\label{F4} \beg{split}& \|\xi\|_{*} \liminf_{k\to\infty}   \W_{p,\ll}(\mu_{k_n}^\theta,\mu^\theta)\ge  \limsup_{k\to\infty} \int_\S  \big|\mu_{k_n}(\e^{\i \|\xi\|_* \theta})- \mu(\e^{\i \|\xi\|_* \theta})\big|\ll(\d\theta)\\
&\ge  \liminf_{k\to\infty} \int_{D_\vv} \Big(\big|\mu_{k_n}(\e^{\i\xi})-\mu (\e^{\i\xi})\big|- \big(\mu+\mu_{k_n}\big)\big(|\e^{\i \|\xi\|_* \theta} -\e^{\i  \xi}|\big)
 \Big)\ll(\d\theta)\\
&\ge \int_{D_\vv}\Big(\vv- \ff \vv 2 \Big) \ll(\d\theta)   =\ff\vv 2\ll(D_\vv)>0,\end{split}\end{equation}
which contradicts to $\W_{p,\ll}(\mu_n,\mu)\to 0$.

(b) Let $\W_{p,\ll}(\mu_n,\mu)\to 0.$ To verify  \eqref{AB0}, let $f_N(r):= (|r|-N)^+,\ r\in \R.$ Then
$$f_N(r)^p\le 2^{p-1} f_N(s)^p+ 2^{p-1} |r-s|^p,\ \ \ r,s\in\R.$$
This implies
\beq\label{SO} \mu^\theta_n(f_N^p)\le 2^{p-1} \mu^\theta(f_N^p) + 2^{p-1} \W_p(\mu_n^\theta,\mu^\theta)^p,\ \ \ n\ge 1,\ \theta\in \S.\end{equation}
Since $\W_{p,\ll}(\mu_n,\mu)\to 0,$ for any $\vv>0$, we find $n_\vv\in \N$ such that
$$\W_{p,\ll}(\mu_n,\mu)^p=\int_{\S} \W_p(\mu_n^\theta,\mu^\theta)^p\ll(\d\theta)\le  \vv,\ \ \ n\ge n_\vv.$$
Combining this with \eqref{SO}, and noting that
$\mu_n\in \scr P_p(\BB)$ implies
$$\lim_{N\to\infty}\int_{\S} \mu_n^\theta(f_N^p)\ll(\d\theta)\le  \lim_{N\to\infty} \mu_n \big([(\|\cdot\|-N)^+]^p\big)=0, $$
we derive
\beg{align*}&\lim_{N\to\infty} \sup_{n\ge 1} \int_{\S} \mu^\theta_n(f_N^p)\ll(\d\theta) \\
&\le 2^{p-1}\lim_{N\to\infty} \sum_{i=1}^{n_\vv} \int_{\S} \mu^\theta_i(f_N^p)\ll(\d\theta) + 2^{p-1}\vv = 2^{p-1}\vv.\end{align*}
By letting $\vv\downarrow 0$ we prove  \eqref{AB0}.

(c) Assuming $\mu_n\fdd\mu$, we prove the equivalence of  \eqref{AB0} and \eqref{AB0'}.
By $\mu_n\fdd\mu$, we have
\beq\label{PU}\lim_{n\to\infty} \mu_n^\theta(|\cdot|^p\land N) =\mu^\theta(|\cdot |^p\land N),\ \ \ \theta\in \S,\ N\in (0,\infty).\end{equation}
Combining this with
\beg{align*}&\bigg|\int_\S \mu_n^\theta(|\cdot|^p)\ll(\d\theta)- \int_\S \mu^\theta(|\cdot|^p)\ll(\d\theta)\bigg|\\
&\le\int_\S \big|\mu_n^\theta(|\cdot|^p\land N) - \mu^\theta(|\cdot|^p\land N)\big|\ll(\d\theta)\\ &\quad + \int_\S\Big[\mu_n^\theta\big((|\cdot|^p-N)^+\big)   +\mu^\theta\big((|\cdot|^p-N)^+\big)\Big]\ll(\d\theta),\end{align*}
and applying the dominated convergence theorem, we obtain
\beg{align*} &\limsup_{n\to\infty} \bigg|\int_\S \mu_n^\theta(|\cdot|^p)\ll(\d\theta)- \int_\S \mu^\theta(|\cdot|^p)\ll(\d\theta)\bigg|\\
&\le  \sup_{n\ge 1}  \int_\S\Big[\mu_n^\theta\big((|\cdot|^p-N)^+\big) +\mu^\theta\big((|\cdot|^p-N)^+\big)\Big]\ll(\d\theta),\ \ N\in (0,\infty).\end{align*}
By letting $N\to\infty,$ we deduce \eqref{AB0'} from \eqref{AB0}.

On the other hand, if  \eqref{AB0'} holds, then by combining with  \eqref{PU}  we obtain
 \beq\label{LM0} \beg{split}&\lim_{n\to\infty} \int_{\S} \mu_n^\theta\big((|\cdot|^p-N)^+\big)\ll(\d\theta) = \lim_{n\to\infty} \int_{\S} \big\{\mu_n^\theta(|\cdot|^p) -\mu_n^\theta
\big((|\cdot|^p\land N) \big) \big\}\ll(\d\theta)\\
&= \int_{\S} \big\{\mu^\theta(|\cdot|^p) -\mu^\theta
\big((|\cdot|^p\land N) \big) \big\}\ll(\d\theta)=\int_{\S}\mu^\theta\big((|\cdot|^p-N)^+\big)\ll(\d\theta),\ \ N\in (0,\infty).\end{split}\end{equation}
Since $\mu\in \scr P_p(\BB)$, for any $\vv>0$ we find $N_\vv\in\mathbb N$ such that
$$\int_{\S}\mu^\theta\big((|\cdot|^p-N)^+\big)\ll(\d\theta)\le \mu\big((\|\cdot\|^p-N)^+\big)\le \vv,\ \ \ N\ge N_\vv.$$
Combining this with \eqref{LM0},  we find $n_\vv\in\mathbb N$ such that
\beq\label{G1} \int_{\S} \mu_n^\theta\big((|\cdot|^p-N_\vv)^+\big)\ll(\d\theta)\le 2\vv,\ \ n\ge n_\vv.\end{equation}
Since $\mu_n\in\scr P_p(\BB)$, we find $N_\vv'\ge N_\vv$ such that
$$\sum_{i=1}^{n_\vv} \int_{\S} \mu_i^\theta\big((|\cdot|^p-N)^+\big)\ll(\d\theta)\le \vv,\ \ \ N\ge N_\vv'.$$
Combining this with \eqref{G1} we derive
\beg{align*} &\sup_{n\ge 1} \int_{\S} \mu_n^\theta\big((|\cdot|^p-N)^+\big)\ll(\d\theta) \\
&\le \sum_{i=1}^{n_\vv} \int_{\S} \mu_i^\theta\big((|\cdot|^p-N)^+\big)\ll(\d\theta) + \sup_{n\ge n_\vv}   \int_{\S} \mu_n^\theta\big((|\cdot|^p-N)^+\big)\ll(\d\theta)
\le 3\vv,\ \ \ N\ge N_\vv'.\end{align*}
Therefore, \eqref{AB0} holds.

(d) Assuming  $\mu_n\fdd\mu$, we show that  \eqref{AB0}  and \eqref{AB0'} imply $\W_{p,\ll}(\mu_n,\mu)\to 0$.

 For  any $N,n\ge 1$ and $\theta\in\S$, let $\pi_{N,n}^\theta\in \C(\mu_n^\theta,\mu^\theta)$ such that
\beg{align*} \W_p^{N,n}(\mu_n^\theta,\mu^\theta)   &:=   \bigg(\inf_{\pi\in  \C(\mu_n^\theta,\mu^\theta) } \int_{\R\times\R} (|r-s|\land N)^p \pi (\d s,\d r)\bigg)^{\ff 1 p}\\
&=\bigg(\int_{\R\times\R} (|r-s|\land N)^p \pi_{N,n}^\theta(\d s,\d r)\bigg)^{\ff 1 p}.\end{align*}
Since   $\mu_n^\theta\to\mu^\theta$ weakly and  $\sup_{n\ge 1} \W_p^{N,n} (\mu_n^\theta,\mu^\theta)   \le N,$ we have
\beq\label{G2} \sup_{n\ge 1} \W_p^{N,n} (\mu_n^\theta,\mu^\theta)   \le N,\ \ \ \lim_{n\to\infty} \W_p^{N,n}(\mu_n^\theta,\mu^\theta)  =0,\ \ \ N\in\mathbb N,\ \theta\in\S.\end{equation}
Moreover,
\beg{align*} &\W_p(\mu_n^\theta,\mu^\theta)^p\le  \int_{\R\times\R} |r-s|^p \pi_{N,n}^\theta (\d s,\d r)\\
&\le \W_p^{N,n}(\mu_n^\theta,\mu^\theta)^p  + \int_{\R\times\R} |r-s|^p (1_{\{|r|>N/2\}}+1_{\{|s|>N/2\}})  \pi_{N,n}^\theta (\d s,\d r)\\
&\le \W_p^{N,n}(\mu_n^\theta,\mu^\theta)^p  + 2^{p-1} \int_{\R\times\R}(|r|^p+|s|^p) (1_{\{|r|>N/2\}}+1_{\{|s|>N/2\}})  \pi_{N,n}^\theta (\d s,\d r).\end{align*}
Since $\pi_{N,n}^\theta\in \C(\mu_n^\theta,\mu^\theta)$, we have
\beg{align*} &  \int_{\R\times\R}(|r|^p+|s|^p) (1_{\{|r|>N/2\}}+1_{\{|s|>N/2\}})  \pi_{N,n}^\theta (\d s,\d r)\\
&= (\mu_n^\theta+\mu^\theta)(|\cdot|^p1_{\{2|\cdot|>N\}}) + \mu_n^\theta(|\cdot|^p)\mu^\theta(\{2|\cdot|>N\})+  \mu^\theta(|\cdot|^p)\mu_n^\theta(\{2|\cdot|>N\}).\end{align*}
By combing these with \eqref{AB0}, \eqref{AB0'} and \eqref{G2}, we obtain
$$  \limsup_{n\to\infty} \int_{\S}  \W_p(\mu_n^\theta,\mu^\theta)^p\ll(\d\theta)=0. $$
 \end{proof}

\beg{proof}[Proof of Theorem $\ref{T1}(2)$]  Let $\{\mu_n\}_{n\ge 1}\subset \scr P_p(\BB)$ be a Cauchy sequence under $\W_{p,\ll}$ satisfying \eqref{BD}, we aim to find $\mu\in \scr P_p(\BB)$ such that
$\W_{p,\ll}(\mu_n,\mu)\to 0.$

(a) Let $\{e_m\}_{m\ge 1}$ be an orthonormal basis of $\BB$. For any $\mu\in \scr P(\BB)$ and non-empty finite subset $T$ of $\mathbb N$ (denoted by $T\Subset \mathbb N$),
we define the marginal distribution $\mu^T$ of $\mu$   as
$$\mu^T:=\mu\circ \Phi_T^{-1},\ \ \ \Phi_T: \BB\to\R^T,\ \Phi_T(x)=(\<x,e_i\>)_{i\in T}.$$
Since $\{\mu_n\}_{n\ge 1}\subset \scr P_p(\BB)$ is a Cauchy sequence under $\W_{p,\ll}$,
  we conclude  that
  $\{\mu_n\big(\e^{\i \<\xi, \cdot\>})\}_{n\ge 1}$ is a Cauchy sequence for any $\xi\in \BB$. Otherwise, there exist $0\ne \xi\in\BB,\vv>0$  and $m_k,n_k\uparrow\infty$ as $k\uparrow\infty$ such that
  $$\big|\mu_{n_k}\big(\e^{\i \<\xi, \cdot\>})-\mu_{m_k}\big(\e^{\i \<\xi, \cdot\>})\big|\ge \vv,\ \ \ k\ge 1.$$
  Noting that $\BB=\BB^*$ as $\BB$ is a Hilbert space, and  that \eqref{BD} implies \eqref{KZ},
  by using $\mu_{m_k}$ replacing $\mu$ in the  argument leading to \eqref{F4}, we derive
  $$\|\xi\|\liminf_{k\to\infty} \W_{p,\ll}(\mu_{n_k},\mu_{m_k})\ge \ff\vv 2 \ll(D_\vv) >0,$$
  which is impossible as $\{\mu_n\}_{n\ge 1}$ is a $\W_{p,\ll}$-Cauchy sequence. This contradiction implies that
$$\varphi(\xi):= \lim_{n\to\infty} \mu_n\big(\e^{\i \<\xi, \cdot\>}),\ \ \ \xi\in \BB$$ exists. By \eqref{BD},
$$\lim_{\xi\to 0} |\varphi(\xi)-1|\le \lim_{\xi\to 0}\sup_{n\ge 1} \|\xi\|\mu_n(\|\cdot\|) =0,$$ so that $\varphi$ is continuous at $\xi=0$.
Then according to  L\'evy's continuity theorem,      for any  $T \Subset  \mathbb N$, there exists a unique probability measure $\mu^{T}\in \scr P(\R^T)$
such that
$\mu_n^{T} \to \mu^{T}\ \text{weakly}.$
By the marginal property
$$\mu_n^{T_1}=\mu_n^{T_2}(\cdot\times \R^{T_2-T_1}),\ \ \ T_1\subset T_2\Subset \mathbb N,$$
the same property holds for $\{\mu^{T}\}_{T\Subset \mathbb N}.$ Therefore, by  the Kolmogorov consistent  theorem, there exists a unique probability
measure $\bar\mu$ on $\R^{\mathbb N}$ equipped with the $\si$-field induced by measurable cylindrical  functions, such that
for any $ T\Subset \mathbb N$, its marginal distribution on $\R^T$ is $\mu^T$.

(b) For any $m\in\mathbb N$, by $\mu_n^{\{1,\cdots, m\}}\to \mu^{\{1,\cdots,m\}}$ weakly and \eqref{BD},   we have
$$\int_{\R^{\mathbb N}} \Big(\sum_{i=1}^m x_i^2 \Big)^{\ff p 2} \bar\mu(\d x)\le  \liminf_{n\to\infty} \int_{\R^{\mathbb N}}\Big(\sum_{i=1}^m x_i^2\Big)^{\ff p 2} \mu_n^{\{1,\cdots,m\}}(\d x)
\le \liminf_{n\to \infty}  \mu_n(\|\cdot\|^p)<\infty.$$
By letting $m\to\infty$ we obtain
$$\int_{\R^{\mathbb N}} \Big(\sum_{i=1}^\infty x_i^2 \Big)^{\ff p 2} \bar\mu(\d x)\le \liminf_{n\to \infty}  \mu_n(\|\cdot\|^p)<\infty.$$
So, $\bar\mu$-a.s.
$$\R^{\mathbb N}\ni x=(x_i)_{i\ge 1}\mapsto \Phi  (x):=\sum_{i=1}^\infty x_i e_i\in \BB,$$
and
$\mu:= \bar\mu\circ \Phi^{-1} \in \scr P_p(\BB)$ satisfies
 $\mu_n\fdd\mu.$   By Theorem \ref{T1}(1), it remains to verify \eqref{AB0}.

Noting that
$$(|r|^p-N)^+\le 2^{p-1} (|s|^p-N)^++ 2^{p-1}|r-s|^p,\ \ \ r,s\in\R,$$ we obtain
\beq\label{*}\mu_n^\theta\big((|\cdot|^p-N)^+\big)\le 2^{p-1} \mu_m^\theta\big((|\cdot|^p-N)^+\big)+ 2^{p-1}\W_p(\mu_n^\theta,\mu_m^\theta)^p,\ \ n\ge m,\ \theta\in\S,\ N\ge 1.\end{equation}
Since $\{\mu_n\}_{n\ge 1}$ is $\W_{p,\ll}$-Cauchy, for any $\vv>0$ we find $m\ge 1$ such that
\beq\label{*2} \sup_{n\ge m} \W_p(\mu_n^\theta,\mu_m^\theta)^p\le 2^{-p}\vv.\end{equation}
Moreover, by $\{\mu_n\}_{n\ge 1}\subset\scr P_p(\BB)$, we find $N_\vv\ge 1$ such that
$$\sup_{1\le i\le m} \int_\S \mu_i^\theta\big((|\cdot|^p-N)^+\big)\ll(\d\theta)\le \sup_{1\le i\le m} \mu_i\big((\|\cdot\|^p-N)^+\big)\le 2^{-p}\vv,\ \ N\ge N_\vv. $$
Combining this with \eqref{*} and \eqref{*2}, we derive
$$\sup_{n\ge 1} \int_\S \mu_i^\theta\big((|\cdot|^p-N)^+\big)\ll(\d\theta)\le   \vv,\ \ \ N\ge N_\vv,$$
so that \eqref{AB0} holds.
  \end{proof}

Finally, to prove Theorem \ref{T1}(3), we need the following lemma.
\beg{lem}\label{L1} Let $\BB=\R^d$ and $\ll\in \scr P(\S)$ with support  containing a  basis $\{e_i\}_{1\le i\le d}$ of  $\R^d$. Then for any $p\in [1,\infty)$  there exists a constant $c\in (0,\infty)$ such that
$$\mu(\|\cdot\|^p)\le c \int_{\S} \mu^\theta(|\cdot|^p)\ll(\d\theta),\ \ \ \mu\in \scr P_p(\R^d).$$
\end{lem}

\beg{proof} Since  $\{e_i\}_{1\le i\le d}$ is a basis of $\R^d$, we find a constant $c_1\in (0,\infty)$
$$|x|^p\le c_1 \Big(\sum_{i=1}^d \<x,e_i\>^2\Big)^{\ff p 2},\ \ x\in\R^d.$$
Noting that
$$2|x|^2 \sum_{i=1}^d |e_i-\theta_i|^2 + 2 \sum_{i=1}^d \<\theta_i,x\>^2\ge 2 \sum_{i=1}^d \big(\<x,e_i-\theta_i\>^2+\<\theta_i,x\>^2\big)\ge \sum_{i=1}^d \<x,e_i\>^2,$$
we find   $c_2\in (0,\infty)$ such that  for any $\{\theta_i\}_{1\le i\le d}\subset \S$,
\beq\label{AF} \beg{split}&|x|^p \le  c_1 \bigg(2|x|^2\sum_{i=1}^d |e_i-\theta_i|^2+ 2\sum_{i=1}^d \<\theta_i,x\>^2\bigg)^{\ff p 2}\\
&\le c_2  |x|^p \sum_{i=1}^d |e_i-\theta_i|^p+ c_2\sum_{i=1}^d|\<\theta_i,x\>|^p,\ \ x\in\R^d.\end{split}\end{equation}
Let
$$ B_i:= \bigg\{\theta\in \S:\ |\theta-e_i|^p\le \ff 1 {2c_2d}\bigg\},\ \ 1\le i\le d.$$
Then \eqref{AF} implies
\beq\label{BJ} |x|^p\le 2c_2 \sum_{i=1}^d |\<\theta_i,x\>|^p,\ \ \ x\in \R^d,\ \theta_i\in B_i,\ 1\le i\le d.\end{equation}
Integrating with respect to $\mu(\d x)\prod_{i=1}^d \ll(\d\theta_i)$ over $\R^d\times \prod_{i=1}^d B_i$ gives
\beg{align*}\mu(|\cdot|^p)\prod_{i=1}^d \ll(B_i) &\le 2 c_2 \sum_{i=1}^d \bigg(\prod_{j\ne i} \ll(B_j)\bigg)\int_{B_i} \mu^{\theta_i}(|\cdot|^p)\ll(\d\theta_i)\\
&\le 2 c_2 \sum_{i=1}^d \bigg(\prod_{j\ne i} \ll(B_j)\bigg)\int_{\S} \mu^{\theta}(|\cdot|^p)\ll(\d\theta).\end{align*}
Since each $e_i$ is included in the support of $\ll$, we have
$\ll(B_i)>0$ for $1\le i\le d.$
Therefore,  the desired estimate holds for
$$c:=    2c_2 \sum_{i=1}^d \ff 1 {\ll(B_i)}\in (0,\infty).$$

\end{proof}

\beg{proof}[Proof of Theorem $\ref{T1}(3)$]  Let $\BB=\R^d$. Since $\W_{p,\ll}\le \W_p$, and $\scr P_p(\R^d)$ is separable under $\W_p$, it is separable under $\W_{p,\ll}$ as well.
Let $\{\mu_n\}_{n\ge 1}\subset \scr P_p(\R^d)$ be a $\W_{p,\ll}$-Cauchy sequence.
By Theorem \ref{T1}(2), to show that $\mu_n$ converges to some $\mu\in \scr P_p(\R^d)$,
we only need to verify \eqref{BD}.  Since  $\{\mu_n\}_{n\ge 1}\subset \scr P_p(\R^d)$ is a $\W_{p,\ll}$-Cauchy sequence, we have
$$\sup_{n\ge 1}  \int_{\S} \mu_n^\theta(|\cdot|^p)\ll(\d\theta)<\infty.$$
This together  Lemma \ref{L1} implies \eqref{BD}. So, $(\scr P_p(\R^d), \W_{p,\ll})$ is a Polish space.

Now, let $\{\mu_n,\mu\}_{n\ge 1}\subset \scr P_p(\R^d)$, we intend to show that   $\W_{p,\ll}(\mu_n,\mu)\to 0$ if and only if $\W_p(\mu_n,\mu)\to 0$.
Since $\W_{p,\ll}\le \W_p$, it suffices to prove that $\W_{p,\ll}(\mu_n,\mu)\to 0$ implies  $\W_p(\mu_n,\mu)\to 0$.
By \eqref{BJ}, there exists a constant $c\in (0,\infty)$ such that for any $N\in (1,\infty)$,
$$(\|x\|^p-N)^+\le c \sum_{i=1}^d \big(|\<x, \theta_i\>|^p-c^{-1} N)^+,\ \ x\in\R^d,\ \theta_i\in B_i.$$
Integrating both sides with respect to $\mu_n(\d x)\times \prod_{i=1}^d 1_{B_i}(\theta_i)\ll(\d\theta_i)$ gives
$$\bigg(\prod_{i=1}^d \ll(B_i)\bigg) \mu_n\big((\|x\|^p-N)^+\big)\le c \sum_{i=1}^d \bigg(\prod_{j\ne i} \ll(B_j)\bigg)
 \int_{B_i} \mu^{\theta_i}(|\cdot|^p-c^{-1}N)^+\ll (\d\theta_i).$$
 Therefore,
 $$\mu_n\big(\|\cdot\|^p-N)^+\big)\le c   \bigg(\sum_{i=1}^d \ff 1 {\ll(B_i)}\bigg) \int_\S\mu_n^\theta\big((|\cdot|^p-c^{-1}N)^+\big) \ll(\d \theta).$$
Since $\W_{p,\ll}(\mu_n,\mu)\to 0$ implies \eqref{BD} as shown above,  we obtain
 \beq\label{LO} \lim_{N\to\infty} \sup_{n\ge 1} \mu_n\big(\|\cdot\|^p-N)^+\big)=0.\end{equation}
 Moreover, by Theorem \ref{T1}(1) for $\BB=\R^d$,
$\W_{p,\ll}(\mu_n,\mu)\to 0$ implies that $\mu_n\to\mu$ weakly, which together with  \eqref{LO} implies  $\W_p(\mu_n,\mu)\to 0$.
  \end{proof}

 The next result provides a H\"older type inequality for the SW distance.

 \beg{prp}\label{PI} For any $\ll\in \scr P_0(\S)$ and $1\le p_1\le p\le p_2<\infty$,
 $$\W_{p,\ll}(\mu,\nu)\le \W_{p_1,\ll}(\mu,\nu)^{\ff{p_1(p_2-p)}{p(p_2-p_1)}} \W_{p_2,\ll}(\mu,\nu)^{\ff{p_2(p-p_1)}{p(p_2-p_1)}},\ \ \ \mu,\nu\in \scr P(\BB).$$
  \end{prp}

 \beg{proof} It suffices to prove for $p_2>p>p_1$, so that
 \beq\label{QD}\vv:= \ff{p_1(p_2-p)}{p_2-p_1}\in (0,p_1),\ \ \ \ff{p_1(p-\vv)}{p_1-\vv}=p_2.\end{equation}
 For any $\theta\in \S$, consider the inverse cumulative distribution functions of $\mu^\theta$:
 $$g_{\mu^\theta}(r):= \inf\big\{s\in\R: \mu^\theta((-\infty,s))\ge r\big\},\ \ \ r\in (0,1),$$
 and define $g_{\nu^\theta}$ in the same way.
According to \cite[Theorem 2]{Ru}  for $\si(x,y)=|x-y|^p$, we have
$$\W_p(\mu^\theta,\nu^\theta)^p= \int_0^1 |g_{\mu^\theta}(r)-g_{\nu^\theta}(r)|^p \d r.$$
Therefore, by H\"older's inequality and \eqref{QD}, we obtain
\beg{align*} &\W_p(\mu^\theta,\nu^\theta)^p =   \int_0^1 |g_{\mu^\theta}(r)-g_{\nu^\theta}(r)|^\vv |g_{\mu^\theta}(r)-g_{\nu^\theta}(r)|^{p-\vv}  \d r \\
&\le \bigg(\int_0^1 |g_{\mu^\theta}(r)-g_{\nu^\theta}(r)|^{p_1}\d r\bigg)^{\ff\vv{p_1}}  \bigg(\int_0^1 |g_{\mu^\theta}(r)-g_{\nu^\theta}(r)|^{\ff{p_1(p-\vv)}{p_1-\vv}}\d r\bigg)^{\ff{p_1-\vv}{p_1}}  \\
&= \W_{p_1}(\mu^\theta,\nu^\theta)^{\ff{p_1(p_2-p)}{p_2-p_1}} \W_{p_2}(\mu^\theta,\nu^\theta)^{\ff{p_2(p-p_1)}{p_2-p_1}}.\end{align*}
Integrating with respect to $\ll(\d\theta)$ and applying H\"older's inequality, we derive
\beg{align*} &\W_{p,\ll}(\mu,\nu)^p\le \int_\S \W_{p_1}(\mu^\theta,\nu^\theta)^{\ff{p_1(p_2-p)}{p_2-p_1}} \W_{p_2}(\mu^\theta,\nu^\theta)^{\ff{p_2(p-p_1)}{p_2-p_1}}\ll(\d \theta)\\
&\le \bigg(\int_\S \W_{p_1}(\mu^\theta,\nu^\theta)^{p_1}\ll(\d\theta)\bigg)^{\ff{p_2-p}{p_2-p_1}}\bigg(\int_\S \W_{p_2}(\mu^\theta,\nu^\theta)^{p_2}\ll(\d\theta)\bigg)^{\ff{p-p_1}{p_2-p_1}}\\
&=  \W_{p_1,\ll}(\mu,\nu)^{\ff{p_1(p_2-p)}{p_2-p_1}} \W_{p_2,\ll}(\mu,\nu)^{\ff{p_2(p-p_1)}{p_2-p_1}}.\end{align*}
This finishes the proof.
 \end{proof}

\section{Application to  exponential ergodic Markov processes}

From now on, we assume that $(X_t)_{t\ge 0}$  is a  Markov process on $\BB$.
For any $\nu\in \scr P(\BB)$, let $\E^\nu$ be the expectation of the Markov process $X_t$ on $\BB$ with initial distribution $\nu$. We   simply denote
$\E^{\dd_x}=\E^x$ for the Dirac measure $ \dd_x$ at $x\in\BB$. Then
\beq\label{AAB}\E^\nu=\int_\BB \E^x \nu(\d x),\ \ \ \nu\in \scr P(\BB).\end{equation}
Let $\B_b(\BB)$ be the space of all bounded measurable functions on $\BB$.
The associated Markov semigroup $(P_t)_{t\ge 0}$ is defined as
$$P_tf(x):= \E^x[f(X_t)],\ \ \ f\in \B_b(\BB),\ t\ge 0,\ x\in\BB,$$
 which  extends uniquely to a $C_0$-contraction semigroup in $L^p(\mu)$ for any $p\in [1,\infty].$

 In this section, we assume the $L^2$-exponential ergodicity.

 \beg{enumerate}\item[$(A_1)$] $P_t$ has an invariant probability measure $\mu$, and there exist $c_0,\kk_0\in (0,\infty)$   such that
$$\|P_t-\mu\|_{L^2(\mu)}:=\sup_{\|f\|_{L^2(\mu)}\le 1} \|P_tf-\mu(f)\|_{L^2(\mu)}\le c_0\e^{-\kk_0 t},\ \ \ t\ge 0.$$
\end{enumerate}

Under $(A_1)$, we study the convergence rate under $\W_{p,\ll}$ for   the empirical measure
$$\mu_t:=\ff 1 t \int_0^t \dd_{X_s}\d s\to\mu\ \text{as}\ t\to\infty,$$
and  the time-discrete empirical measure
$$\tt\mu_n :=\ff 1 n\sum_{i=1}^n \dd_{X_i}\to \mu \ \text{as}\ n\to\infty.$$

We first explain that in may cases   the study can be reduced  to the stationary case where $\mu$ is the initial distribution.
For any initial distribution $\nu\le c \mu$ for some constant $c>0$,  \eqref{AAB} implies $\E^\nu \le c  \E^\mu.$
 Moreover,   if $P_t$ has heat kernel $p_t(x,\cdot)$ with respect to $\mu$, then
\beg{align*}&\E^x[\W_{p,\ll}(\mu_{t+1},\mu)]= \int_\BB p_1(x,y) \E^y [\W_{p,\ll}(\mu_t,\mu)]\mu(\d y)\\
&\le  \bigg(\int_\BB p_1(x,y)^2\mu(\d y)\bigg)^{\ff 1 2}\ss{\E^\mu [\W_{p,\ll}(\mu_t,\mu)^2]},\ \ x\in\BB, t>0.\end{align*}
In the following, we present another result which reduces   arbitrary initial distributions to the stationary case, without using the heat kernel.

\beg{prp}\label{PR} Let $X_t$ be the Markov process on $\BB$ with stationary distribution $\mu$, and let $p\in [1,\infty)$. If for some constants $c,\aa\in (0,\infty)$ and any $x,y\in \BB$, there exists a coupling
$(X_t^x,X_t^y)$ of the Markov processes with $X_0^x=x$ and $X_0^y=y$ such that
\beq\label{PR0}\E\|X_t^x-X_t^y\|^p\le c\e^{-\aa t} |x-y|^p, \ \ \ t\ge 0,\end{equation}
then for any $\ll\in \scr P_0(\BB)$,
\beq\label{PR1} \E^\nu\big[\W_{p,\ll}(\mu_t,\mu)^p\big]\le 2^{p-1}\Big(\ff{c} {\aa t} \W_p(\nu,\mu)^p+ \E^\mu \big[\W_{p,\ll}(\mu_t,\mu)^p\big]\Big),\ \ \ t>0,\end{equation}
\beq\label{PR2} \E^\nu\big[\W_{p,\ll}(\tt\mu_n,\mu)^p\big]\le 2^{p-1}\Big(\ff{c} {\aa n} \W_p(\nu,\mu)^p+ \E^\mu \big[\W_{p,\ll}(\tt\mu_n,\mu)^p\big]\Big),\ \ \ n\in \mathbb N,\end{equation}

\end{prp}

\beg{proof} We only prove \eqref{PR1} as the other estimate can be proved in the same way.
Let
$$\mu_t^x=\ff 1 t \int_0^t \dd_{X_s^x}\d s,\ \ \ \mu_t^y=\ff 1 t \int_0^t \dd_{X_s^y}\d s,\ \ \ t>0.$$ Then
$$\E\big[\W_{p,\ll}(\mu_t^x,\mu_t^y)^p\big]\le \ff 1 t \int_0^t \E|X_s^x-X_s^y|^p\d s\le \ff{2^{p-1} c}{\aa t} |x-y|^p,\ \ \ x,y\in\BB,\ t>0.$$ Combining this with the triangle inequality, we derive
\beg{align*}&\E^x\big[\W_{p,\ll}(\mu_t,\mu)^p\big]=\E\big[\W_{p,\ll}(\mu_t^x,\mu)^p\big] \le 2^{p-1}\E \Big(\W_{p,\ll}(\mu_t^y,\mu)^p+ \W_{p,\ll}(\mu_t^x,\mu_t^y)^p\Big)\\
&\le 2^{p-1}\E^y \big[\W_{p,\ll}(\mu_t,\mu)^p\big]+    \ff{2^{p-1} c}{\aa} |x-y|^p,\ \ \  x,y\in\BB,\ t>0.\end{align*}
Then \eqref{PR1} follows by integrating with respect to the optimal coupling $\pi\in \C(\nu,\mu)$ with
$$\W_p(\mu,\nu)^p= \int_{\BB\times \BB}\|x-y\|^p\pi(\d x,\d y).$$
 \end{proof}


 \subsection{Dimension-free convergence rate in $\W_{p,\ll}$}
 For any $\theta\in\S$ and $t>0$, consider the cumulative distribution function  of $\mu^\theta$:
 $$F^\theta(r):= \mu^\theta\big((-\infty,r)\big),\ \  \ \ r\in\R.$$

 \beg{thm}\label{W1} Let  $X_t$ be a Markov process satisfying $(A_1)$. Then for any $p\in [1,\infty)$ and $\ll\in \scr P_0(\BB),  $
 \beq\label{QA1} \beg{split}&\E^\mu\big[\W_{p,\ll}(\mu_t,\mu)^{2p}\big]\le \int_\S \E^\mu\big[\W_p(\mu_t^\theta,\mu^\theta)^{2p}\big]\ll(\d\theta)\\
&\le \ff{c_0p^2 2^{2p-1}}{\kk_0 t} \int_\S\bigg(\int_\R |r|^{p-1} \ss{F^\theta(r)(1-F^\theta(r))}\, \d r\bigg)^2\ll(\d\theta),\ \ t>0,\end{split}\end{equation}
and for $K_p:= p^2 4^{p-1}+  c_0p^2 2^{2p-1} \sum_{i=1}^\infty \e^{-\kk_0 i},$
\beq\label{QA2}\beg{split}& \E^\mu\big[\W_{p,\ll}(\tt\mu_n,\mu)^{2p}\big]\le \int_\S \E^\mu\big[\W_p(\tt\mu_n^\theta,\mu^\theta)^{2p}\big]\ll(\d\theta)\\
&\le \ff {K_p}n   \int_\S\bigg(\int_\R |r|^{p-1} \ss{F^\theta(r)(1-F^\theta(r))}\, \d r\bigg)^2\ll(\d\theta),\ \ n\in\mathbb N.\end{split} \end{equation}
\end{thm}

 \beg{proof} Let $F_t^\theta$ be the  cumulative distribution function  of $\mu_t^\theta$. Then
\beq\label{FT}F_t^\theta(r)= \ff 1 t\int_0^t 1_{(-\infty,r)}(\theta(X_s))\d s,\ \ \ r\in\R.\end{equation}
 By \cite[Proposition 7.4]{BL},  we have
\beg{align*}&\W_p(\mu_t^\theta,\mu^\theta)^p\le p 2^{p-1} \int_\R |r|^{p-1} \big|F_t^\theta (r)-F^\theta (r)\big|\d r\\
&= p 2^{p-1} \int_\R |r|^{p-1} \bigg|\ff 1 t\int_0^t \big[1_{(-\infty,r)}(\theta(X_s))-F^\theta (r)\big]\d s\bigg|\d r.\end{align*}
So,   for any positive measurable function $h$ on $\R$, by Schwarz's inequality we derive
 \beq\label{KT}  \beg{split} \W_p(\mu_t^\theta,\mu^\theta)^{2 p} \le  &p^2 4^{p-1} \bigg(\int_\R |r|^{2(p-1)}h(r)\d r\bigg)\\
 &\times \int_\R \ff 1 {h(r)} \bigg|\ff 1 t\int_0^t \big[1_{(-\infty,r)}(\theta(X_s))-F^\theta(r)\big] \d s\bigg|^2\d r.\end{split}\end{equation}
It is easy to see that the function
$$ g_r(x):= 1_{(-\infty,r)}(\theta(x))-F^\theta(r),\ \ \ x\in\BB$$
 satisfies
 \beq\label{GR}\beg{split}&\mu(g_r)= \mu\big(\{x\in\BB: \theta(x)<r\}\big)-F^\theta(r)= 0,\\
 &\mu(g_r^2)= F^\theta(r)- F^\theta(r)^2= F^\theta(r)\big(1-F^\theta(r)\big),\ \ \ r\in\R.\end{split}\end{equation}
By  the Markov property of $X_s$ and that    $\mu$ is the invariant probability measure of $X_s$,    we obtain
 \beq\label{KT2}\beg{split}&  \E^\mu \bigg(\ff 1 t\int_0^t \big[1_{(-\infty,r)}(\theta(X_s))-F^\theta(r)\big] \d s\bigg)^2 \\
 & = \ff 2 {t^2} \int_0^t \d t_1 \int_{t_1}^t \E^\mu\big[g_r(X_{t_1}) g_r(X_{t_2})\big]\d t_2\\
 &= \ff 2 {t^2} \int_0^t \d t_1 \int_{t_1}^t \E^\mu\big[g_r(X_{t_1}) \big(P_{t_2-t_1} g_r\big)(X_{t_1})\big]\d t_2\\
 &\le  \ff 2 {t^2} \int_0^t \d t_1 \int_{t_1}^t \| g_r\|_{L^2(\mu)} \big\|P_{t_2-t_1} g_r\big\|_{L^2(\mu)}\d t_2.\end{split}\end{equation}
Combining this with   \eqref{KT}, \eqref{GR} and  $(A_1)$, we arrive at
 $$\E^\mu\big[\W_p(\mu_t^\theta,\mu^\theta)^{2p}\big] \le \ff{c_0p^2 2^{2p-1}}{\kk_0 t} \bigg(\int_\R |r|^{2(p-1)}h(r)\d r\bigg)\int_\R \ff {F^\theta(r)(1-F^\theta(r))} {h(r)}\d r.$$
 Taking $h(r)= |r|^{1-p}\ss{F^\theta(r)(1-F^\theta(r))}$ gives
$$\E^\mu\big[\W_p(\mu_t^\theta,\mu^\theta)^{2p}\big] \le \ff{c_0p^2 2^{2p-1}}{\kk_0 t} \bigg(\int_\R |r|^{p-1} \ss{F^\theta(r)(1-F^\theta(r))}\,\d r\bigg)^2.$$
Thus, \eqref{QA1} holds.

Next, the  cumulative distribution function  of $\tt\mu_n^\theta$ is
 $$\tt\mu_n^\theta((-\infty,r)) =\ff 1 n\sum_{i=1}^n   1_{(-\infty,r)}(\theta(X_i)),\ \ \ r\in\R. $$   Repeating the above argument for $\tt\mu_n$ in place of $\mu_t$ leads to
\beg{align*} &p^{-2} 4^{1-p} \E^\mu\big[\W_p(\tt\mu_n^\theta,\mu^\theta)^{2p}\big]\\
&\le \bigg(\int_\R |r|^{2(p-1)} h(r)\d r\bigg)\E^\mu\int_\R \ff 1 {h(r)} \bigg(\ff 1 n \sum_{i=1}^n g_r(X_i)  \bigg)^2\d r\\
&=   \bigg(\int_\R |r|^{2(p-1)} h(r)\d r\bigg) \bigg[\int_\R \ff{ \|g_r\|_{L^2(\mu)}^2 }{n h(r)} \d r+ \ff 2 {n^2}  \sum_{1\le i<j\le n} \int_\R\ff{c_0 \|g_r\|_{L^2(\mu)}^2\e^{-(j-i)\kk_0}}{h(r)}\d r\bigg]\\
&\le \bigg(\ff 1 n + \ff {2c_0}{n} \sum_{i=1}^\infty \e^{-\kk_0 i} \bigg) \bigg(\int_\R |r|^{p-1} \ss{F^\theta(r)(1-F^\theta(r))} \d r\bigg)^2.\end{align*}
Therefore, \eqref{QA2} holds.
 \end{proof}

 \beg{cor}\label{C1} Let $p\in [1,\infty)$ and $M_{q}(\mu)= \mu(\|\cdot\|^{q})<\infty$ for some $q>2p$, such that
 $$C_{p,q}(\mu):= \Big( M_1(\mu) + \ff {M_{  q}(\mu)}{ q}\Big) \int_0^\infty \ff {2r^{2(p-1)}} {1+r^{q-1}}\d r<\infty.$$
 Then for any Markov process satisfying $(A_1)$,
  \beq\label{NM0} \beg{split}  &\sup_{\ll\in\scr P_0(\S)}\E^\mu\big[\W_{p,\ll}(\mu_t,\mu)^{2p}\big]\le   \ff{2c_0C_{p,q}(\mu)}{\kk_0 t},\ \  \ \ t>0,\\
  &\sup_{\ll\in\scr P_0(\S)} \E^\mu\big[\W_{p,\ll}(\tt\mu_n,\mu)^{2p}\big]\le   \ff {K_0 C_{p,q}(\mu)}n,\ \ \ \ n\in\mathbb N.\end{split}\end{equation}
 \end{cor}
 \beg{proof}
 Noting that for any $q\ge 1$ we have
\beq\label{NN0} \beg{split}&\mu(\|\cdot\|^q)\ge \int_\R |r|^q \mu^\theta(\d r) = q \int_0^\infty \big[1- F^\theta(r)\big] r^{q-1}\d r +  q \int_{-\infty}^0  F^\theta(r)  |r|^{q-1}\d r\\
&\ge q \int_\R F^\theta(r)\big(1-F^\theta(r)\big) |r|^{q-1}\d r,\end{split}\end{equation}
 by Schwarz's inequality we obtain
\beq\label{C1*}\beg{split}& \bigg(\int_\R |r|^{p-1} \ss{F^\theta(r)(1-F^\theta(r))} \d r\bigg)^2\\
&\le  \bigg(\int_\R F^\theta(r)(1-F^\theta(r)) (1+|r|^{q-1}) \d r\bigg) \int_\R \ff {|r|^{2(p-1)}} {1+|r|^{q-1}} \d r\\
& \le C_{p,q}(\mu).\end{split}\end{equation}
Then \eqref{NM0} follows from   Theorem \ref{W1}.
 \end{proof}

 \paragraph{Remark 3.1.}
By Kantorovich's dual formula and the central limit theorem due to \cite{Wu},
\beg{align*} &\liminf_{t\to\infty}  t\,\E^\mu \big[\W_1(\mu_t^\theta, \mu^\theta)^2\big]= \liminf_{t\to \infty}  t\, \E^\mu \sup_{\|f'\|_\infty\le 1}\bigg|\ff 1 t \int_0^t \big[(f\circ\theta)(X_s)-\mu^\theta(f)\big]\d s\bigg|^2\\
&\ge \liminf_{t\to \infty} \sup_{\|f'\|_\infty\le 1} t\, \E^\mu\bigg|\ff 1 t \int_0^t \big[(f\circ\theta)(X_s)-\mu(f\circ \theta)\big]\d s\bigg|^2 \\
&= 2\sup_{\|f'\|_\infty\le 1, \mu^\theta(f)=0}\int_0^\infty \mu  \big( ( f\circ\theta)   P_s (f\circ\theta) \big)\d s >0,\ \ \ \theta\in\S. \end{align*}
Then, for any $\ll\in \scr P_0(\S)$,
\beq\label{ADS} \liminf_{t\to\infty}  t \E^\mu \big[\W_{1,\ll}(\mu_t, \mu)^2\big]>0.\end{equation}
 The same holds for $ n \E^\mu \big[\W_{1,\ll}(\tt\mu_n, \mu)^2\big] $ as $n\to\infty$. So,
 by \eqref{NM0} for $p=1$, if $M_q(\mu)<\infty$ for some $q>2$, then
    for any Markov process $X_t$ on $\BB$ satisfying $(A_1)$, we have the sharp asymptotic  formulas
$$\E^\mu\big[\W_{1,\ll} (\mu_t,\mu)^2\big] \sim t^{-1},\ \ \ \E^\mu\big[\W_{1,\ll} (\tt\mu_n,\mu)^2\big] \sim n^{-1}.$$
The estimate \eqref{NM0} for $p>1$ is dimension-free but  less sharp.

 In  the next part, we intend to  derive sharp convergence rate  for $ \E^\mu[\W_{2,\ll}(\mu_t,\mu)^2]$
and $ \E^\mu[\W_{2,\ll}(\tt\mu_n,\mu)^2]$ as $t,n\to\infty$.

 \subsection{Sharp dimension-free convergence rate in $\W_{2,\ll}$}

 To estimate $\E^\mu\big[\W_{2,\ll} (\mu_t,\mu)^2\big],$ we  assume that for any $\theta\in \S$,   there exists an open interval $I_\theta$ such that $\mu^\theta(I_\theta)=1$ and
 $$\rr_\theta :=\ff{\d\mu^\theta}{\d r} \in C^1\big(I_\theta; (0,\infty)\big).$$  Let $X_t^\theta$ be the diffusion process generated by the operator $L^\theta$:
$$L^\theta f(r) := f''(r)+ \ff{\rr_\theta'(r) }{\rr_\theta(r)} f'(r),\ \ \ r\in I_\theta, $$
  with reflecting boundary if $I_\theta\ne\R$  and the boundary can be reached by the diffusion.  Define the associated diffusion semigroup
$$P_t^\theta f(r):=\E\big[f(X_t^\theta)|X_0^\theta=r\big],\ \ \ \ t\ge 0,\ f\in \B_b(I_\theta),\ r\in I_\theta.$$
Since $L^\theta$ is symmetric in $L^2(\mu^\theta),$ $(P_t^\theta)_{t\ge 0}$ extends uniquely to a symmetric diffusion semigroup on $L^2(\mu^\theta)$.
We assume that $P_t^\theta$ has a heat kernel $p_t^\theta(r,s)$ with respect to $\mu^\theta$ such that
$$P_t^\theta f(r)= \int_{I_\theta} p_t^\theta(r,s) f(s)\mu^\theta(\d s),\ \ \ t>0,\ f\in\B_b(I_\theta),\ r\in I_\theta.$$
  Besides $(A_1)$, we make the following assumption on $\mu^\theta$ and $P_t^\theta$, which is modified form \cite{W25}.

  \beg{enumerate}\item[$(A_2)$] There exist  measurable functions
  $$\aa, \ k: \ \S\to (0,\infty),\ \ \ \  \gg: \S\times (0,\infty)\to (0,\infty),$$ where $\gg(\theta,t)$ is decreasing in $t$, such that for any $\theta\in \S$,
\beq\label{PK} \mu^\theta(f^2)\le \ff 1 { \aa(\theta) }\mu^\theta(|f'|^2),\ \ \ \ f\in C_b^1(I_\theta),\ \mu^\theta(f)=0,\end{equation}
  \beq\label{PK2}  \int_{I_\theta} P_t^\theta|\cdot-r|^2(r) \mu^\theta(\d r)\le k(\theta) t,\ \ \ t\in (0,1],\end{equation}
  \beq\label{PK3} \int_{I_\theta} p_{t}^\theta(r,r)\mu^\theta(\d r) \le \gg(\theta,t),\ \ \ t\in (0,\infty).\end{equation}
  \end{enumerate}
These conditions refer to the Poincar\'e inequality and  short time behaviors of   one-dimensional diffusions, and thus  can be easily verified in applications, see Propositions \ref{PP1}-\ref{PP3} below.
Under $(A_1)$ and $(A_2)$  we have the following estimates, where the orders in  \eqref{SP} are sharp according to \eqref{ADS}.

\beg{thm}\label{T3} Assume $(A_1)$ and $(A_2)$, let
\beg{align*}&H(\theta,\vv):= \int_0^\infty s^{-\ff 1 2} \e^{-\ff{s}2\aa(\theta)} \ss{\gg(\theta, 2\vv+s)}\d s,\\
& \ K_0:= 1+ 2c_0\sum_{i=1}^\infty \e^{-\kk_0 i}. \end{align*}
Then  $K_0\le  1+\ff{2c_0}{\kk_0},$ and
\beq\label{SP0} \E^\mu \big[\W_{2,\ll}(\mu_t,\mu)^2\big]\\
 \le \int_\S\inf_{\vv \in (0,1]} \bigg(\ss{k(\theta) \vv}+ \ff{  \ss{8c_0} H(\theta,\vv)} {\GG(\ff 1 2)  \ss{\kk_0 t}} \bigg)^2  \ll(\d \theta), \ \ t\in (0,\infty), \end{equation}
\beq\label{SP02} \E^\mu \big[\W_{2,\ll}(\tt\mu_n,\mu)^2\big]
 \le    \int_\S\inf_{\vv \in (0,1]} \bigg(\ss{k(\theta) \vv}+ \ff {2\ss{ K_0 } H(\theta,\vv)} {\GG(\ff 1 2) \ss n}  \bigg)^2\ll(\d \theta),\ \ n\in\mathbb N.
 \end{equation}
Consequently, if $C(\ll):=\int_\S H(\theta,0)^2\ll(\d\theta)<\infty$, then for any $t>0$ and $n\in\mathbb N,$
\beq\label{SP}\E^\mu \big[\W_{2,\ll}(\mu_t,\mu)^2\big]\le  \ff {8 c_0C(\ll)}{\GG(\ff 1 2)^2 \kk_0 t},\ \
 \E^\mu \big[\W_{2,\ll}(\tt\mu_n,\mu)^2\big]\le   \ff {4K_0C(\ll)}{n\GG(\ff 1 2 )^2}.\end{equation}
  \end{thm}

\beg{proof}   We only need to prove \eqref{SP0} and \eqref{SP02}, which imply \eqref{SP}  with $\vv \downarrow 0$.

(a) By \cite[Theorem 2]{Ledoux} and \eqref{PK}, for any probability $\nu\in \scr P(I_\theta)$ with density $f:=\ff{\d\nu}{\d\mu^\theta}\in L^2(\mu^\theta)$,
we have
\beq\label{WW0}\W_2(\nu,\mu^\theta)^2\le 4\big\|(-L^\theta)^{-\ff 1 2 } (f-1)\big\|^2_{L^2(\mu^\theta)} =
  \bigg\|\ff 2 {\GG(\ff 1 2)}\int_0^\infty s^{-\ff 1 2} P_{s}^\theta (f-1)\d s\bigg\|_{L^2(\mu^\theta)}^2.\end{equation}
Since the projected empirical measure
$$\mu_t^\theta= \ff 1 t \int_0^t \dd_{\theta(X_s)}\d s$$
  is singular with respect to $\mu^\theta$, we make the following regularization by using the heat kernel $p_\vv^\theta$ for $\vv\in (0,1]$:
\beq\label{WW1} \mu_{t,\vv}^\theta(\d r):= f_{t,\vv}(r) \mu^\theta(\d r),\ \ \ f_{t,\vv}(r):= \ff 1 t \int_0^t p_\vv^\theta\big(\theta(X_s),r\big)\d s,\ \ \ t>0,\ r\in I_\theta.\end{equation}
By the triangle inequality, we obtain
\beq\label{TR} \W_2(\mu_t^\theta,\mu^\theta)^2\le (1+\vv_0^{-1}) \W_2(\mu_{t,\vv}^\theta,\mu^\theta)^2+ (1+\vv_0) \W_2(\mu_t^\theta, \mu_{t,\vv}^\theta)^2,\ \ \vv,\vv_0\in (0,1].\end{equation}
Below we estimate these two terms   respectively.

(b) Noting that $\mu^\theta(p_\vv^\theta(r,\cdot)^2)=p_{2\vv}^\theta(r,r)$,
by Jensen's inequality  we obtain
\beq\label{WW3}\beg{split} &\E^\mu\big[\mu^\theta (f_{t,\vv}^2)\big]\le \ff 1 t \E^\mu \int_0^t \mu^\theta\Big(p_\vv^\theta\big(\theta(X_s),\cdot\big)^2 \Big)\d s\\
&\le \ff 1 t  \int_0^t  \E^\mu\big[p_{2\vv}^\theta\big(\theta(X_s),\theta(X_s)\big)\big]\d s,\ \ \ \theta\in\S, \ \vv>0.\end{split}\end{equation}
Since $\mu^\theta=\mu\circ \theta^{-1}$ and $\mu$ is the invariant probability measure of $X_s$, we have
$$\E^\mu\big[p_{2\vv}^\theta\big(\theta(X_s),\theta(X_s)\big)\big]=\int_\BB p_{2\vv}^\theta\big(\theta(x),\theta(x)\big)\mu(\d x)= \int_{I_\theta} p_{2\vv}^\theta(r,r)\mu^\theta(\d r).$$
Combining this with \eqref{WW3} and \eqref{PK3} we derive  $\E^\mu\big[\mu^\theta (f_{t,\vv}^2)\big]<\infty$, so that   \eqref{WW0} yields
\beq\label{X}\E^\mu\big[\W_2(\mu_{t,\vv}^\theta,\mu^\theta)^2\big] \le   \bigg(\ff 2 {\GG(\ff 1 2)} \int_0^\infty s^{-\ff 1 2} \big\|P_s^\theta(f_{t,\vv}-1)\big\|_{L^2(\mu^\theta)}\d s\bigg)^2.\end{equation}
By the definition of $f_{t,\vv}$ and noting that $p_s^\theta$ is the symmetric kernel of $P_s^\theta$ with respect to $\mu^\theta$, we obtain
$$P_s^\theta(f_{t,\vv}-1)= P_{\ff{s}2}^\theta (f_{t,\ff s2+\vv}-1).$$ This and \eqref{X} imply
\beq\label{X1} \E^\mu\big[\W_2(\mu_{t,\vv}^\theta,\mu^\theta)^2\big]  \le   \bigg(\ff 2 {\GG(\ff 1 2)}\int_0^\infty s^{-\ff 1 2}
\big\|P_{\ff{s}2}^\theta(f_{t,\ff{s}2+\vv}-1)\big\|_{L^2(\mu^\theta)}\d s\bigg)^2.\end{equation}
By   the Poincar\'e inequality \eqref{PK}, we obtain
$$\big\|P_{\ff{s}2}^\theta(f_{t,\ff s2+\vv}-1)\big\|_{L^2(\mu^\theta)}\le \e^{-\ff{s}2\aa(\theta)} \|f_{t,\ff s2+\vv}-1\|_{L^2(\mu^\theta)},$$
which together with   \eqref{X1} and Schwarz's inequality gives
 \beq\label{WF0}\beg{split} &\W_2(\mu_{t,\vv}^\theta,\mu^\theta)^2
\le \bigg(\ff 2 {\GG(\ff 1 2)}\int_0^\infty s^{-\ff 1 2} \e^{-\ff s 2\aa(\theta)} \big\|f_{t,\ff{s}2+\vv}-1\big\|_{L^2(\mu^\theta)}\d s\bigg)^2\\
&\le \ff 4 {\GG(\ff 1 2)^2}\bigg(\int_0^\infty s^{-\ff 1 2} \e^{-\ff s 2\aa(\theta)} \ss{\gg(\theta,2\vv+s)}\,\d s\bigg) \int_0^\infty \ff{s^{-\ff 1 2} \e^{-\ff s 2\aa(\theta)}}{\ss{\gg(\theta,2\vv+s)}} \big\|f_{t,\ff s2+\vv}-1\big\|_{L^2(\mu^\theta)}^2\d s.\end{split} \end{equation}
 Moreover, by   \eqref{WW1}  we obtain
 $$\big\|f_{t,\ff s2+\vv}-1\big\|_{L^2(\mu^\theta)}^2=\ff 2 {t^2} \int_0^t \d t_1 \int_{t_1}^t \d t_2\int_{I_\theta} \big(p_{\ff{2\vv+s}2}^\theta(\theta(X_{t_1}),r)-1\big)\big(p_{\ff{2\vv+s}2}^\theta(\theta(X_{t_2}),r)-1\big)\mu^\theta(\d r),$$
 so that
 \beq\label{WF1}\beg{split} &\E^\mu\Big[\big\|f_{t,\ff s2+\vv}-1\big\|_{L^2(\mu^\theta)}^2\Big]\\
& =\ff 2 {t^2} \int_0^t \d t_1 \int_{t_1}^t \d t_2\int_{I_\theta} \E^\mu\Big[\big(p_{\ff{2\vv+s}2}^\theta(\theta(X_{t_1}),r)-1\big)\big(p_{\ff{2\vv+s}2}^\theta(\theta(X_{t_2}),r)-1\big)\Big]\mu^\theta(\d r).
 \end{split}\end{equation}
By $\mu^\theta=\mu\circ\theta^{-1}$, $\mu$ is the invariant probability measure of $X_t$, and $p_t^\vv$ is the symmetric heat kernel of $P_t^\theta$ with respect to $\mu^\theta$, we derive
\beq\label{YPP1} \beg{split}&\E^\mu \Big[\big|p_{\ff{2\vv+s}2}^\theta(\theta(X_t),r)-1\big|^2\Big]= \big\|p_{\ff{2\vv+s}2}^\theta(\theta(\cdot),r)-1\big\|_{L^2(\mu)}^2\\
& = \big\|p_{\ff{2\vv+s}2}^\theta(\cdot,r)-1\big\|_{L^2(\mu^\theta)}^2 = p_{s+2\vv}(r,r)-1,\ \ t\ge 0.\end{split}\end{equation}
Combining this with the Markov property of $X_t$,     Schwarz's inequality and $(A_1)$, we arrive at
\beq\label{YPP2} \beg{split}&\E^\mu \Big[\big(p_{\ff{2\vv+s}2}^\theta(\theta(X_{t_1}),r)-1\big)\big(p_{\ff{2\vv+s}2}^\theta(\theta(X_{t_2}),r)-1\big)\Big]\\
&=\E^\mu \Big[\big(p_{\ff{2\vv+s}2}^\theta(\theta(X_{t_1}),r)-1\big) P_{t_2-t_1} \big(p_{\ff{2\vv+s}2}^\theta(\theta(\cdot),r)-1\big)(X_{t_1})\Big]\\
&\le   \big\|p_{\ff{2\vv+s}2}^\theta(\cdot,r)-1 \big\|_{L^2(\mu^\theta)} \big\|P_{t_2-t_1}  (p_{\ff{2\vv+s}2}^\theta(\theta(\cdot),r)-1\big)\big\|_{L^2(\mu)}\\
&\le c_0\e^{-\kk_0(t_2-t_1)} \big\|p_{\ff s2+\vv}(\theta(\cdot),r)-1\big\|_{L^2(\mu^\theta)}^2\\
&= c_0\e^{-\kk_0(t_2-t_1)} \big(p_{s+2\vv}(r,r)-1\big),\ \ \ t_2\ge t_1\ge 0. \end{split}\end{equation}
This together with \eqref{WF0}  and \eqref{WF1} yields
\beq\label{WF2} \beg{split} \E^\mu\Big[\W_2\big(\mu_{t,\vv}^\theta,\mu^\theta\big)^2\Big] &\le \ff {8c_0}{\kk_0\GG(\ff 1 2)^2 t}  \bigg(\int_0^\infty s^{-\ff 1 2} \e^{-\ff s 2\aa(\theta)} \ss{\gg(\theta,2\vv+s)}\,\d s\bigg)^2\\
&= \ff {8c_0H(\theta,\vv)^2}{\kk_0\GG(\ff 1 2)^2 t}.\end{split}\end{equation}

(c)  It is easy to see that for any $t,\vv>0$,
$$\pi_{t,\vv}(\d r_1, \d r_2):= \ff 1 t \int_0^t \big(\dd_{\theta(X_s)}(\d r_1) p_\vv^\theta(\theta(X_s), r_2)\mu^\theta(\d r_2)\big)\d s\in \C(\mu_t^\theta,\mu_{t,\vv}^\theta).$$
So,
\beg{align*} & \W_2(\mu_t^\theta,\mu_{t,\vv}^\theta)^2\le \int_{I_\theta\times I_\theta} |r_1-r_2|^2 \pi_{t,\vv}(\d r_1,\d r_2)\\
&= \ff 1 t \int_0^t |\theta(X_s)-r|^2 p_\vv^\theta(\theta(X_s), r) \mu^\theta(\d r) =\ff 1 t \int_0^t P_\vv^\theta|\theta(X_s)-\cdot|^2(\theta(X_s))\d s.\end{align*}
Since with initial distribution $\mu$, the law of $\theta(X_s)$ is $\mu^\theta$, this together with \eqref{PK2} implies
\beq\label{YPP3} \E^\mu\big[\W_2(\mu_t^\theta,\mu_{t,\vv}^\theta)^2\big]\le k(\theta)\vv,\ \ \ \theta\in \S,\ \vv\in (0,1].\end{equation}
By combining this with \eqref{TR} and \eqref{WF2},  we obtain
\beg{align*} \E^\mu\big[\W_2(\mu_t^\theta,\mu^\theta)^2\big] &\le \inf_{\vv\in (0,1)}\inf_{\vv_0\in (0,1)} \Big\{\big(1+\vv_0^{-1}\big)\cdot \ff {8c_0H(\theta,\vv)^2}{\kk_0\GG(\ff 1 2)^2 t}
+\big(1+\vv_0\big)k(\theta) \vv\Big\}\\
&= \inf_{\vv\in (0,1)} \bigg(\ff{H(\theta,\vv)\ss{8c_0} }{\GG(\ff 1 2) \ss{\kk_0 t}} +\ss{k(\theta)\vv}\bigg)^2.\end{align*}
Therefore,  \eqref{SP0} holds.

  The proof of \eqref{SP02} is similar to that of \eqref{SP0}. Let
  $$\tt f_{n,s}:=\ff 1 n\sum_{i=1}^n p_s^\theta(\theta(X_i), r),\ \ \ \ s> 0,\ n\in\mathbb N.$$
  Then the same reason leading to \eqref{WF0} implies
  \beq\label{WF02} \beg{split}\W_2(\tt\mu_{n,\vv}^\theta,\mu^\theta)^2
  \le \,&\ff 4 {\GG(\ff 1 2)^2}\bigg(\int_0^\infty s^{-\ff 1 2} \e^{-\ff s 2\aa(\theta)} \ss{\gg(\theta,2\vv+s)}\,\d s\bigg) \\
  &\quad\times \int_0^\infty \ff{s^{-\ff 1 2} \e^{-\ff s 2\aa(\theta)}}{\ss{\gg(\theta,2\vv+s)}}
  \|\tt f_{n,\ff {s+2\vv}2}-1\|_{L^2(\mu^\theta)}^2\d s.
  \end{split}\end{equation}
  Noting that
 \beq\label{DP} \beg{split} &\E^\mu\Big[\big\|\tt f_{n,\ff{2\vv +s}2}-1\big\|_{L^2(\mu^\theta)}^2\Big]\\
 & =\ff 1 {n^2} \sum_{i=1}^n \int_{I_\theta} \E^\mu\Big[\big(p_{\ff{2\vv+s}2}^\theta(\theta(X_i),r)-1\big)^2\Big]\mu^\theta(\d r)\\
  &\qquad  +\ff 2 {n^2} \sum_{i=1}^n \sum_{j=i+1}^n  \big(p_{\ff{2\vv+s}2}^\theta(\theta(X_i),r)-1\big) \big(p_{\ff{2\vv+s}2}^\theta(\theta(X_j),r)-1\big)\mu^\theta(\d r),\end{split}\end{equation}
  by \eqref{YPP1} and   \eqref{YPP2}    we derive
  \beg{align*} &\E^\mu\Big[\big\|\tt f_{n,\ff{2\vv +s}2}-1\big\|_{L^2(\mu^\theta)}^2\Big] \le \int_{I_\theta} \Big[
  \ff 1 {n}  \big(p_{ 2\vv+s }^\theta(r,r)-1\big)  +
  \ff 2 {n }   \sum_{j= 1}^\infty  c_0 \e^{-\kk_0 j}   \big(p_{2\vv+s}^\theta(r ,r)-1\big) \Big]\mu^\theta(\d r)\\
  &\le  \ff {K_0\gg (\theta,2\vv+s)} n.\end{align*}
  Combining this with \eqref{WF02},   we arrive at
\beq\label{GM}  \E^\mu \big[\W_{2}(\tt\mu_{n,\vv}^\theta,\mu^\theta)^2\big]
 \le     \ff {4  K_0   H(\theta,\vv)^2} {\GG(\ff 1 2)^2  n}.\end{equation}
 Similarly to     \eqref{YPP3}, we have
 \beq\label{YPP3N} \E^\mu\big[\W_2(\tt\mu_n, \tt\mu_{n,\vv})^2\big]\le k(\theta)\vv,\ \ \ n\in\mathbb N,\ \vv\in (0,1].\end{equation}
 This together with \eqref{GM} and the triangle inequality implies   \eqref{SP02}.
\end{proof}

To apply Theorem \ref{T3},  we present below some explicit conditions ensuring \eqref{PK}-\eqref{PK3}, where the first result
is due to the Hardy inequality in \cite[Theorem 3.1]{Chen00}, which implies $\aa(\theta)>0$ if $\rr_\theta(r)\sim \e^{-c |r|^k}$ for some constants $c>0$ and $k\ge 1$.

\beg{prp}\label{PP1}    For any $\theta\in \S$, let $  I_\theta =(a_\theta, b_\theta)$.   Then   for any $r_0\in I_\theta$,
$$4\aa(\theta)\ge  \min\bigg\{\inf_{r\in (r_0,b_\theta)} \ff 1 {(\int_{r_0}^{r} \rr_\theta(s)\d s )\int_r^{b_\theta} \rr_\theta(s)\d s},\ \inf_{r\in (a_\theta,r_0)} \ff 1 {(\int_r^{r_0} \rr_\theta(s)\d s )\int_{a_\theta}^r \rr_\theta(s)\d s}\bigg\}.$$
\end{prp}
\beg{proof} For any $r\in I_\theta$, let
\beg{align*}&\aa_+(r):=\inf\bigg\{\int_{r}^{b_\theta} f'(s)^2 \rr_\theta(s)\d s:\ f\in C_b^1([r,b_\theta)),\ \int_{r}^{b_\theta} f(s)\rr_\theta(s)\d s=0\bigg\},\\
&\aa_-(r):=\inf\bigg\{\int_{a_\theta}^r f'(s)^2 \rr_\theta(s)\d s:\ f\in C_b^1((a_\theta,r]),\ \int_{a_\theta}^r f(s)\rr_\theta(s)\d s=0\bigg\},\end{align*}
Then $\aa_+(r)$ is increasing in $r$,  $\aa_-(r)$ is decreasing in $r$, and there exists a unique $r_\theta\in I_\theta$ such that
$$\aa(\theta)= \aa_+(r_\theta)=\aa_-(r_\theta).$$
By the Hardy inequality in \cite[Theorem 3.1]{Chen00},  we have
$$ \aa_+(r_\theta)\ge \ff 1 4 \inf_{r\in (r_\theta,b_\theta)} \ff 1 {(\int_{r_\theta}^{r} \rr_\theta(s)\d s )(\int_r^{b_\theta} \rr_\theta(s)\d s)},$$
$$ \aa_-(r_\theta)\ge \ff 1 4 \inf_{r\in (a_\theta,r_\theta)} \ff 1 {(\int_{r}^{r_\theta} \rr_\theta(s)\d s )(\int_{a_\theta}^{r} \rr_\theta(s)\d s)}.$$
Then the proof is finished.
\end{proof}

\beg{prp}\label{PP2} In general,  $\eqref{PK2}$ holds for
$$ k(\theta)=4 + 2 \big\|\rr_\theta'/\rr_\theta\big\|_{L^2(\mu^\theta)}^2.$$
If $I_\theta$ is bounded, $\log \rr_\theta \in C^2(I_\theta)$ and there exists a constant $K_\theta\in \R$ such that
\beq\label{CV}  \ff{\d^2}{\d r^2} \log \rr_\theta(r)  \le K_\theta,\ \ \ r\in I_\theta,\end{equation}
then $\eqref{PK2}$ holds for
$$k(\theta)= \e^{K_\theta^+}\bigg(2+ |I_\theta|\int_{I_\theta} |\rr_\theta'|(r) \d r\bigg).$$
\end{prp}

\beg{proof}
Let $X_t^\theta$ be the diffusion process generated by $L^\theta$, we may find a Brownian motion on $\R$ such that
$$X_t^\theta-X_0^\theta= \ss 2\, W_t+ \int_0^t \ff{\rr_\theta'}{\rr_\theta}(X_s^\theta)\d s,\ \ \ t\ge 0.$$ So, letting $X_0^\theta$ have initial distribution $\mu^\theta$ which is the invariant probability measure of $X_s^\theta$, we derive
\beg{align*} &\int_\R P_t^\theta|\cdot-r|^2(r) \mu^\theta(\d r)= \E |X_t^\theta-X_0^\theta|^2 \\
&\le 4 \E[W_t^2] + 2 \E \bigg( \int_0^t \ff{\rr_\theta'}{\rr_\theta}(X_s^\theta)\d s\bigg)^2\\
&\le 4t + 2 t \int_0^t \E\Big|\ff{\rr_\theta'}{\rr_\theta}(X_s^\theta)\Big|^2\d s
= 4 t + 2t \big\|\rr_\theta'/\rr_\theta\big\|_{L^2(\mu^\theta)}^2,\ \ \ t\in (0,1]. \end{align*}
So, $\eqref{PK2}$ holds for $k(\theta)=4 + 2 \big\|\rr_\theta'/\rr_\theta\big\|_{L^2(\mu^\theta)}^2.$

Next, let $I_\theta$ be bounded and \eqref{CV} hold. Then for   any $r,s\in I_\theta, s\ne r$, we have
\beg{align*} &L^\theta |\cdot-r|(s) =    {\rm sgn}(s-r) \Big(\ff{\d}{\d s} \log \rr_\theta(s)\Big) \\
&= {\rm sgn}(s-r) \Big(\ff{\d}{\d r} \log \rr_\theta(r)\Big) + {\rm sgn}(s-r) \int_r^s \Big(\ff{\d^2}{\d t^2} \log \rr_\theta(t)\Big)\d t\\
&\le  \Big|\ff{\d}{\d r} \log \rr_\theta(r)\Big|+ K_\theta |s-r|.\end{align*}
Hence,
$$L^\theta |\cdot-r|^2(s) \le 2 + |I_\theta|  \Big|\ff{\d}{\d r} \log \rr_\theta(r)\Big|+ K_\theta^+  |s-r|^2,\ \ s\in I_\theta,$$
so that
$$P_t^\theta |\cdot-r|^2(r)\le t \e^{K_\theta^+} \bigg(2+  |I_\theta|  \ff{|\rr_\theta'|}{\rr_\theta}(r) \bigg),\ \ \ r\in I_\theta,\ t\in [0,1].$$
This implies $\eqref{PK2}$  for $k(\theta)= \e^{K_\theta^+}\big(2+ |I_\theta|\int_{I_\theta} |\rr_\theta'|(r) \d r\big).$ \end{proof}

Finally, we present a concrete choice of $\gg(\theta,t)$ in the condition \eqref{PK3}.

\beg{prp}\label{PP3} If there exists $K_\theta\in\R$ such that $\eqref{CV}$ holds, then
\beq\label{HT} p_{t}^\theta(r,r)\le \ff 1 {\int_{I_\theta}  \exp [-\ff{K_\theta|r-s|^2}{1-\exp[-K_\theta t]} ]\mu^\theta(\d s)},\ \ t>0,\end{equation}
so that $\eqref{PK3}$ holds for
$$\gg(\theta,t):= \int_{I_\theta}\ff{\mu^\theta(\d r)}{\int_{I_\theta} \exp\big(-\ff{K_\theta|r-s|^2}{1-\exp[-K_\theta t]}\big)\mu^\theta(\d s)},\ \ \ t>0,$$
where $\ff{K_\theta|r-s|^2}{1-\exp[-K_\theta t]}:= \ff {|r-s|^2}t$ if $K_\theta=0.$
\end{prp}

\beg{proof} By \cite[Lemma 2.1]{W97} for $P_t=P_t^\theta, K=K_\theta$,  $g(s)=\e^{-K_\theta s}$ and $f(r):= p_t^\theta(r,\cdot)$, we obtain
$$p_{2t}^\theta(r,r)^2= \big(P_t^\theta f(r)\big)^2\le \big(P_t^\theta f^2 (s)\big)\exp\Big[\ff{K_\theta|r-s|^2}{1-\exp[-2K_\theta t]}\Big],
\ \ \ \ r,s\in\R,\ t>0.$$
Since $\mu^\theta$ is the invariant probability measure of $P_t^\theta$, we obtain
$$p_{2t}^\theta(r,r)^2\int_{I_\theta}  \exp\Big[-\ff{K_\theta|r-s|^2}{1-\exp[-2K_\theta t]}\Big]\mu^\theta(\d s)\le \int_\R P_t^\theta f^2 (s)\mu^\theta(\d s)
= \mu^\theta(f^2)= p_{2t}^\theta(r,r).$$
Then
$$p_{2t}^\theta(r,r)\le \ff 1 {\int_{I_\theta}  \exp\Big[-\ff{K_\theta|r-s|^2}{1-\exp[-2K_\theta t]}\Big]\mu^\theta(\d s)},\ \ \ t>0,\ r\in I_\theta.$$
This is equivalent to \eqref{HT}.
\end{proof}

To conclude this subsection, we present two examples to  illustrate Theorem \ref{T3} by verifying assumption $(A_2)$  using the above propositions.
We first consider   exponential ergodic Markov processes on $\BB$ with a Gaussian measure as invariant probability measure. In particular, when $\BB=\R^d$ and $X_t$ is the Ornstein-Uhlenbeck process,
 the dimension-free convergence rate of  $\E[\W_{2,\ll}(\mu_t,\mu)^2]$   reaches  the sharp order  $t^{-1}$   uniformly in $\ll\in \scr P_0(\S)$.

\beg{exa} Let $\mu$ be a Gaussian measure on $\BB$ with mean $x_0:=\int_\BB x\mu(\d x)\in\BB$, such that
$$\si_\theta:=\bigg(\int_\BB\big|\theta(x)- \theta(x_0)\big|^2\mu(\d x)\bigg)^{\ff 1 2}= \bigg(\int_\R \big|r- \theta(x_0)\big|^2 \mu^\theta(\d r)\bigg)^{\ff 1 2}$$
is bounded in $\theta\in \S$. Let $\si:=\sup_{\theta\in\S} \si_\theta.$ Then the following assertions hold.
\beg{enumerate} \item[$(1)$] There exists a universal constant $c\in (0,\infty)$ such that for any Markov process $X_t$ on $\BB$ satisfying $(A_1)$,
\beg{align*}&\sup_{\ll\in \scr P_0(\S)}\E^\mu\big[\W_{2,\ll}(\mu_t,\mu)^2\big] \le c \Big(1+c_0^{\ff 1 2}\kk_0^{-\ff 1 2}\Big) \Big( \ff{1+ \si^3}t+\ff{\si[\log (1+t)]^2}t\Big),\ \ \ t>0,\\
&\sup_{\ll\in \scr P_0(\S)}\E^\mu\big[\W_{2,\ll}(\tt\mu_n,\mu)^2\big] \le c \Big(1+c_0^{\ff 1 2}\kk_0^{-\ff 1 2}\Big) \Big( \ff{1+ \si^3}n+\ff{\si[\log (1+n)]^2}n\Big),\ \ \ n\ge 1.\end{align*}
\item[$(2)$] Let $\si\in (0,\infty)$ and let   $X_t$ be  the Ornstein-Uhlenbeck process generated by the operator
$$ L^{\rm ou} f(x):= \DD f(x)- \si^{-2} \<x, \nn f(x)\>,\ \ \ x\in\R^d.$$  Then $\mu$ is the Normal distribution $N(0, \si^{-2} I_d)$  where $I_d$ is the $d\times d$-identity matrix, and
\beg{align*}    \sup_{\ll\in \scr P_0(\S)}\E^\mu\big[\W_{2,\ll}(\mu_t,\mu)^2\big]\le \ff {8\si^2}t  \sum_{i=1}^\infty \ff 1 {i^2},\ \ t>0.  \end{align*}
Moreover, there exists a universal constant $c\in (0,\infty)$ such that
$$ \sup_{\ll\in \scr P_0(\S)}\E^\mu\big[\W_{2,\ll}(\tt\mu_n,\mu)^2\big]\le  \ff {c (1+ \log [n(4\si^2+ 1)])}n,\ \ \ n\in \mathbb N.$$
 \end{enumerate}
 \end{exa}

\beg{proof}  Without loss of generality, we assume that $x_0:=\int_\BB x\mu(\d x)=0$. Then
\beq\label{RR0} \rr_\theta(r):=\ff{\mu^\theta(\d r)}{\d r}= \ff {\exp[-\ff{r^2}{2\si_\theta^2}]} {\ss{2\pi}\,\si_\theta },\ \ \ r\in\R,\ \theta\in\S.\end{equation}
It is well-known that   the   Ornstein-Uhlenbeck operator
$$L^\theta f(r):= f''(r)- \ff{ r f'(r)}{\si_\theta^2},\ \ \ \ r\in\R   $$ has spectral gap   $ \si_\theta^{-2}$, so \eqref{PK} holds for
\beq\label{RR1} \aa(\theta)= \si_\theta^{-2}>0.\end{equation}
Next, by \eqref{RR0}  and Proposition \ref{PP2}, \eqref{PK2} holds for
\beq\label{RR2} k(\theta):= 4+ 2\big\|\rr_\theta'/\rr_\theta\big\|_{L^2(\mu^\theta)}^2=4+ \si_\theta^{-2}\in (0,\infty).\end{equation}
Finally,   noting that \eqref{RR0} implies
$$\ff{\d^2}{\d r^2} \log \rr_\theta(r)  = -\si_\theta^{-2}\le 0,$$
by Proposition \ref{PP3} we derive \eqref{PK3} for
$$ \gg(\theta,t)= \int_\R \ff{\mu^\theta(\d r)}{\int_\R \exp[-\ff{|r-s|^2}t]\mu^\theta(\d s)}=\int_\R\ff{\exp[-\ff{r^2}{2\si_\theta^2}]\d r}{\int_\R \exp[-\ff{|r-s|^2} t- \ff{s^2}{2\si_\theta^2}]\d s}.$$
Noting that
$$-\ff{|r-s|^2} t- \ff{s^2}{2\si_\theta^2}= -\Big(t^{-1}+\ff 1{2} \si_\theta^{-2}\Big) \Big(s-\ff r{1+\ff 1 2 t\si_\theta^{-2}}\Big)^2 -\ff{r^2}{2\si_\theta^2+t},$$
we obtain
\beg{align*} & \gg(\theta,t)=\Big(\ff{t^{-1}+ \ff 1 2\si_\theta^{-2}}{2\pi}\Big)^{\ff 1 2} \int_\R \exp\Big[- r^2\Big(\ff 1 {2\si_\theta^2} - \ff 1 {2\si_\theta^2+t}\Big)\Big]\d r\\
&= \Big(t^{-1}+\ff 1 2\si_\theta^{-2}\Big)^{\ff 1 2}\ss{2\si_\theta^2(1+2\si_\theta^2 t^{-1})}=1+2\si_\theta^2t^{-1}.\end{align*}
Combining this with \eqref{RR1}, we derive
\beg{align*} &H(\theta,\vv)= \int_0^\infty s^{-\ff 1 2} \e^{-\ff{s^2}{2\si_\theta^2}} \ss{1+2\si_\theta^2(2\vv+s)^{-1}}\,\d s\\
&\le \int_0^\vv s^{-\ff 1 2}\big(1+2\si_\theta \vv^{-\ff 1 2}\big)\d s +\int_\vv^1 \big(s^{-\ff 1 2}+ \ss 2 \si_\theta s^{-1}\big)\d s +\ss{1+2\si_\theta^2} \int_1^\infty  \e^{-\ff s{2\si_\theta^2}}\d s\\
&\le 2+ 2  \si_\theta + 2\si_\theta^2 \ss{1+2\si_\theta^2} +\ss 2 \si_\theta \log \vv^{-1},\ \ \vv\in (0,1).\end{align*}
By combining this with \eqref{RR2}, and taking
$$\vv= \ff 1 {(1+t)(4+\si_\theta^{-2})}$$ in  \eqref{SP0} and \eqref{SP02},
  we prove the desired estimates for some universal constant $c\in (0,\infty)$.

(2) Let $\mu=N(0, \si^2 I_d)$. Then   $\mu^\theta= N(0,\si^{-2})$ for any $\theta\in \S:=\{\theta\in\R^d: |\theta|=1\}$,
and $X_t^\theta:=\theta(X_t)$ is the one-dimensional Ornstein-Uhlenbeck process generated by
$$ L_\si:= \ff{\d^2}{\d r^2} -\ff r {\si^2}  \ff{\d}{\d r}\  \text{on}\  \R.$$
It is well-known that   $-L_\si$ is self-adjoint in $L^2(N(0,\si^{-2})$, and its all eigenvalues
$$\aa_i= \si^{-2} i,\ \ \ i\ge 0$$
are simple with unitary eigenfunctions
$$\phi_i(r):=  u_i(\si^{-2}r),\ \ \ i\ge 0,\ r\in \R,$$  where $\{u_i\}_{i\ge 0}$ are the normalized Hermit polynomials with
$$   \ff 1 {\ss{2\pi}} \int_\R u_i(r)^2\e^{-\ff{r^2}2}\d r=1,\ \ i\ge 0.$$
  So, $\{\phi_i\}_{i\ge 0}$ consists an orthonormal basis of $L^2(\R,N(0,\si^2))$, and
$$p_\vv^\theta(r,s)= 1+\sum_{i=1}^\infty \e^{-\si^{-2} i \vv} \phi_i(r)\phi_i(s),\ \ \ \vv>0,\ r,s\in\R.$$
Noting that for any $\theta\in\S$, $X_t^\theta:=\theta(X_t)=\<X_t,\theta\>$ is the diffusion process generated by $L_\si$, and
$$f_{t,\vv}^\theta:=\ff{\d\mu_{t,\vv}^\theta}{\d\mu^\theta}= \ff 1 t\int_0^t p_\vv(X_s^\theta,\cdot)\d s,$$
we obtain
$$4\big\|(-L^\theta)^{-1} (f_{t,\vv}^\theta-1)\big\|_{L^2(\mu^\theta)}^2= \ff 8 {t^2} \sum_{i=1}^\infty \ff {\e^{-2\si^{-2} i\vv}}  {\aa i} \int_0^t \d t_1\int_{t_1}^t \phi_i(\theta(X_{t_1}))  \phi_i(\theta(X_{t_2})) \d t_2.$$
By the Markov property of $\theta(X_t)$ and that $\phi_i$ is the unit eigenfunction of $L_\si$ with eigenvalue $-\si^{-2} i$, we obtain
$$\E^\mu\big[\phi_i(\theta(X_{t_1}))  \phi_i(\theta(X_{t_2}))\big]=\e^{-\si^{-2} i(t_2-t_1)},$$
so that
$$\E^\mu\big[4\big\|(-L^\theta)^{-1} (f_{t,\vv}^\theta-1)\big\|_{L^2(\mu^\theta)}^2\big]\le \ff {8\si^2} {t}\sum_{i=1}^\infty \ff 1 {i^2}.$$
This together with \eqref{WW0} for $\nu=\mu_{t,\vv}^\theta$ implies
$$\E^\mu\big[\W_2(\mu_{t,\vv}^\theta,\mu^\theta)^2\big]\le \ff {8\si^2} {t}\sum_{i=1}^\infty \ff 1 {i^2}.$$
Therefore, by Fatou's lemma,
$$\E^\mu\big[\W_2(\mu_{t}^\theta,\mu^\theta)^2\big] = \E^\mu\Big[\lim_{\vv\to 0} \W_2(\mu_{t,\vv}^\theta,\mu^\theta)^2\Big]\le \ff {8\si^2} { t}\sum_{i=1}^\infty \ff 1 {i^2}.$$

Noting that $\d \tt\mu_{n,\vv}^\theta= \tt f_{n,\vv}\d\mu^\theta$ with
$$\tt f_{n,\vv}^\theta= \ff 1 n \sum_{k=1}^n p_\vv^\theta(X_k^\theta,\cdot),$$
we find a universal constant $c_1\in (0,\infty)$ such that
\beg{align*} &\E^\mu\big[\W_2(\tt\mu_{n,\vv}^\theta,\mu^\theta)^2\big]\le \E^\mu\big[4\big\|(-L^\theta)^{-1} (\tt f_{n,\vv}^\theta-1)\big\|_{L^2(\mu^\theta)}^2\big]\\
&= \sum_{i=1}^\infty \ff{2\e^{-2\si^{-2}i}}{i n^2} \E^\mu \Big(\sum_{k=1}^n \phi_i(\theta(X_k))^2 + \sum_{k=1}^{n-1} \sum_{j=k+1}^n  \phi_i(\theta(X_k))\phi_i(\theta(X_j))\Big)\\
&\le \ff 8 n \sum_{i=1}^\infty  \ff{2\e^{-2\si^{-2}i}}{i n^2}\Big(1+ \sum_{l=1}^\infty \e^{-l i}\Big)\le c_1 \log(1+\si^2\vv^{-1}),\ \ \ \vv\in (0,1].\end{align*}
Combining this with \eqref{YPP3N} for $k(\theta)= 4+\si^{-2}$ due to \eqref{RR2}, by the triangle inequality and taking
$$\vv=\ff 1 {n (4+\si^{-2})},$$
we find a constant $c\in (0,\infty)$ such that
$$ \E^\mu\big[\W_2(\tt \mu_{n}^\theta,\mu^\theta)^2\big]\le  \ff {c (1+ \log [n(4\si^2+ 1)])}n,\ \ \  n\in \mathbb N.$$
Then the proof is finished.
\end{proof}

In the following example, we consider the uniform distribution $\mu$ on the unit ball in $\BB=  \R^d$, and derive sharp dimension-free convergence rate for $\E^\mu[\W_{2,\ll}(\mu_t,\mu)^2]$ and $\E^\mu[\W_{2,\ll}(\tt\mu_n,\mu)^2]$  for any exponentially ergodic Markov processes satisfying $(A_1)$.

\beg{exa} Let $d\ge 2$ and $\mu(\d x)=1_{B_d(1)} \oo_d^{-1} \d x$, where $B_d(1):=\{x\in \R^d: |x|<1\}$  and $\oo_d$ is the volume of $B_d(1)$. Then there exists  $c\in (0,\infty)$ uniformly in $d\ge 2$ such that
for any Markov process satisfying $(A_1)$,
\beg{align*}&   \sup_{\ll\in \scr P_0(\S)}\E^\mu\big[\W_{2,\ll}(\mu_t,\mu)^2\big]\le \ff{c c_0}{ \kk_0 \, t},\ \ t>0,\\
& \sup_{\ll\in \scr P_0(\S)}\E^\mu\big[\W_{2,\ll}(\mu_t,\mu)^2\big]\le  c \Big(1+\ff{c_0}{ \kk_0}\Big) \ff 1n,\ \ n\in\mathbb N.  \end{align*}
  \end{exa}
\beg{proof} It is clear that for any $\theta\in\S$,
\beq\label{AA0} \rr_\theta(r)= \rr(r):= \ff 1 {\oo_d} \int_{\{z\in\R^{d-1}\}} \d z= \ff{\oo_{d-1}(1-r^2)^{\ff{d-1}2}}{\oo_d},\ \ \ r\in I_\theta=(-1,1).\end{equation}
Then
\beq\label{AA}  \ff{\d^2 }{\d r^2}\log \rr(r)= -\ff{d-1}{1-r^2}-\ff{2(d-1)r^2}{(1-r^2)}\le -(d-1),\ \ \ r\in I_\theta.\end{equation}
By e.g. \cite[Proposition 4.8.1]{BGL},  this implies  \eqref{PK} for
\beq\label{AA2} \aa(\theta)\ge d-1.\end{equation}
Moreover, according to  the second assertion in Proposition \ref{PP2},   by \eqref{AA}, $I_\theta=(-1,1)$ and \eqref{AA0}, we derive   \eqref{PK2}   for
$$k(\theta)= 2 +2 \int_{-1}^1 |\rr_\theta'(r)|\d r =
 2 + \ff{2(d-1)\oo_{d-1}}{\oo_d} \int_{-1}^1(1-r^2)^{\ff{d-3}2}|r|\d r<\infty.  $$
Finally, by \eqref{AA} and Proposition \ref{PP3}, \eqref{PK3} holds for
\beq\label{AA4} \gg(\theta,t)= \int_{-1}^1 \ff{(1-r^2)^{\ff{d-1}2} }{\int_{-1}^1 \e^{-\ff{|s-r|^2}t}(1-s^2)^{\ff{d-1}2}\d s}\, \d r,\ \ \ t>0.\end{equation}
If $r\in [\ff 1 4,1)$, then $1-s^2\ge 1-r^2$ for $s\in [r-\ff 1 2 \ss{t\land 1}, r]$, so that
$$\int_{-1}^1 \e^{-\ff{|s-r|^2}t}(1-s^2)^{\ff{d-1}2}\d s\ge \e^{-\ff 1 4} \int_{r- \ff 1 2 \ss{t\land 1}}^r (1-s^2)^{\ff {d-1}2}\d s \ge \ff 1 2 \e^{-\ff 1 4} \ss{t\land 1} (1-r^2)^{\ff{d-1}2}.$$
By the symmetry, the same estimate holds for $r\in (-1,-\ff 14]$, so   that
\beq\label{DD} \int_{-1}^1 \e^{-\ff{|s-r|^2}t}(1-s^2)^{\ff{d-1}2}\d s\ge  \ff 1 2 \e^{-\ff 1 4} \ss{t\land 1} (1-r^2)^{\ff{d-1}2},\ \ \ff 1 4\le |r|<1.\end{equation}
Now, let $r\in [0,\ff 1 4)$ and $r_0= \ss{t}\land \ff 2{d+4}.$ If $s\in [r,r+r_0]$, then
$$\ff{1-s^2}{1-r^2}=1-\ff{s^2-r^2}{1-r^2}\ge 1-\ff{ r_0/2+r_0^2}{15/16}\ge 1- \ff{4}{d+4},$$
so that there exists a universal constant $c_0\in (0,\infty)$ such that
\beg{align*} &\int_{-1}^1 \e^{-\ff{|s-r|^2}t}(1-s^2)^{\ff{d-1}2}\d s\ge \e^{-1}\int_{r}^{r+r_0} (1-s^2)^{\ff{d-1}2}\d s\\
&\ge \e^{-1} r_0 (1-r^2)^{\ff{d-1}2}  \Big(1- \ff{4}{d+4}\Big)^{\ff{d-1}2}\ge c_0 (1-r^2)^{\ff{d-1}2}  \Big(\ss t\land \ff{2}{d+4}\Big).\end{align*}
By the symmetry, the same estimate holds for $r\in (-\ff 1 4,0].$ Combining this with \eqref{AA4} and \eqref{DD}, we find a universal constant $c_1\in (0,\infty)$ such that
\beq\label{CD1}  \gg(\theta,t)\le c_1\big( t^{-\ff 1 2}+d\big),\ \ \ t>0.\end{equation}
Thus, there exists a universal constant $c_2\in (0,\infty)$ such that
$$H(\theta,0)= \int_0^\infty  s^{-\ff 1 2} \e^{-\ff{(d-1)s}2} \ss{c_1\big(s^{-\ff 1 2}+ d\big)}\,\d s\le c_2.$$
Then the desired estimates follow from Theorem \ref{T3}.
\end{proof}

\section{An extension to  ergodic Markov processes}

In this part, we study ergodic Markov processes with convergence rate
\beq\label{XI} (0,1]\ni \xi(t):= \|P_t-\mu\|_{L^\infty(\mu)\to L^2(\mu)}\downarrow 0\ \text{as}\ t\uparrow \infty,\end{equation}
which may be slower than exponential.  This type convergence rate has been characterized  in \cite{RW01}                   by the   weak Poincar\'e inequality,
where a number of examples are presented.

The following result extends Theorem \ref{W1} to Markov processes satisfying \eqref{XI}.

 \beg{thm}\label{W1'} Let  $p\in [1,\infty)$ and  $\ll\in \scr P_0(\BB)$. Then the following estimates hold for  any Markov process $X_t$ satisfying $\eqref{XI}$:
\beq\label{QA1'} \beg{split}&\E^\mu\big[\W_{p,\ll}(\mu_t,\mu)^{2p}\big]\le \int_\S \E^\mu\big[\W_p(\mu_t^\theta,\mu^\theta)^{2p}\big]\ll(\d\theta)\\
&\le \ff{p^2 2^{2p}}{t} \bigg(\int_0^t \xi(s)\d s \bigg)\int_\S\bigg(\int_\R |r|^{p-1} \big[F^\theta(r)(1-F^\theta(r))\big]^{\ff 1 4}\, \d r\bigg)^2\ll(\d\theta),\ \ t>0,\end{split}\end{equation}
\beq\label{QA2'}\beg{split}& \E^\mu\big[\W_{p,\ll}(\tt\mu_n,\mu)^{2p}\big]\le \int_\S \E^\mu\big[\W_p(\tt\mu_n^\theta,\mu^\theta)^{2p}\big]\ll(\d\theta)\\
&\le \bigg(\ff 1 n +\ff 2 {n} \sum_{i=1}^n \xi(i)\bigg)  \int_\S\bigg(\int_\R |r|^{p-1} \big[F^\theta(r)(1-F^\theta(r))\big]^{\ff 1 4}\, \d r\bigg)^2\ll(\d\theta),\ \ n\in\mathbb N.\end{split} \end{equation}
\end{thm}

\beg{proof} By \eqref{KT2} and \eqref{XI}, we obtain
\beg{align*}  &  \E^\mu \bigg(\ff 1 t\int_0^t \big[1_{(-\infty,r)}(\theta(X_s))-F^\theta(r)\big] \d s\bigg)^2 \\
 &\le  \ff {2\| g_r\|_{L^2(\mu)}} {t^2} \int_0^t \d t_1 \int_{t_1}^t  \xi(t_2-t_1)\d t_2\le \ff {2\| g_r\|_{L^2(\mu)}} {t}\int_{0}^t  \xi(s)\d s. \end{align*}
Combining this with   \eqref{KT} and \eqref{GR},    we arrive at
 $$\E^\mu\big[\W_p(\mu_t^\theta,\mu^\theta)^{2p}\big] \le \bigg(\ff{2p^2 2^{2p-1}}{ t} \int_0^t \xi(s)\d s\bigg)\bigg(\int_\R |r|^{2(p-1)}h(r)\d r\bigg)\int_\R \ff {\ss{F^\theta(r)(1-F^\theta(r))}} {h(r)}\d r.$$
 Taking $h(r)= |r|^{1-p} [F^\theta(r)(1-F^\theta(r))]^{\ff 1 4}$ gives
$$\E^\mu\big[\W_p(\mu_t^\theta,\mu^\theta)^{2p}\big] \le \bigg(\ff{2p^2 2^{2p-1}}{ t} \int_0^t \xi(s)\d s\bigg)\bigg(\int_\R |r|^{p-1}\big[F^\theta(r)(1-F^\theta(r))\big]^{\ff 1 4}\,\d r\bigg)^2.$$
Hnece, \eqref{QA1'} holds.

Similarly,
\beg{align*} &p^{-2} 4^{1-p} \E^\mu\big[\W_p(\tt\mu_n^\theta,\mu^\theta)^{2p}\big]\\
&\le \bigg(\int_\R |r|^{2(p-1)} h(r)\d r\bigg)\E^\mu\int_\R \ff 1 {h(r)} \bigg(\ff 1 n \sum_{i=1}^n g_r(X_i)  \bigg)^2\d r\\
&\le   \bigg(\int_\R |r|^{2(p-1)} h(r)\d r\bigg) \bigg[\int_\R \ff{ \|g_r\|_{L^2(\mu)}  }{n h(r)} \d r+ \ff 2 {n^2}  \sum_{1\le i<j\le n} \int_\R\ff{ \|g_r\|_{L^2(\mu)}\xi(j-i)}{h(r)}\d r\bigg]\\
&\le \bigg(\ff 1 n + \ff {2}{n} \sum_{i=1}^n \xi(i) \bigg) \bigg(\int_\R |r|^{p-1} \big[F^\theta(r)(1-F^\theta(r))\big]^{\ff 1 4} \d r\bigg)^2.\end{align*}
Therefore, \eqref{QA2'} holds.
\end{proof}

By \eqref{C1*}, we have the following consequence of Theorem \ref{W1'}.

 \beg{cor}\label{C1'} Let $p\in [1,\infty)$ and $M_{q}(\mu)= \mu(\|\cdot\|^{q})<\infty$ for some $q>4p$, such that
 $$\tt C_{p,q}(\mu):= \Big( M_1(\mu) + \ff {M_{q}(\mu)}{q}\Big)^{\ff 1 2}  \bigg(\int_0^\infty \ff {2 r^{\ff 4 3}}{(1+r^{q-1})^{\ff 1 3}}\d r\bigg)^{\ff 3 2}<\infty.$$
 Then for any Markov process satisfying $\eqref{XI}$, the following estimates hold:
\beg{align*}&  \E^\mu\big[\W_{p,\ll}(\mu_t,\mu)^{2p}\big]\le   \ff{2\tt C_{p,q}(\mu)}{t}\int_0^t \xi(s)\d s,\ \  \ \ t>0,\\
 &  \E^\mu\big[\W_{p,\ll}(\tt\mu_n,\mu)^{2p}\big]\le  \tt C_{p,q}(\mu)\bigg(\ff 1 n +\ff 2 {n} \sum_{i=1}^n \xi(i)\bigg),\ \ \ \ n\in\mathbb N.\end{align*}
 \end{cor}
\beg{proof} By \eqref{NN0} and H\"older's inequality,  we obtain
\beg{align*} &\bigg(\int_\R |r|^{p-1}\big[F^\theta(r)(1-F^\theta(r))\big]^{\ff 1 4}\d r\bigg)^2\\
&\le \bigg(\int_\R (1+ |r|^{q-1}) F^\theta(r)(1-F^\theta(r)) \d r\bigg)^{\ff 12} \bigg(\int_\R \ff{|r|^{\ff 4 3(p-1)}} {(1+ |r|^{q-1})^{\ff 1 3}}\d r   \d r\bigg)^{\ff 3 4}\\
& \le \tt C_{p,q}(\mu).\end{align*} Then the proof is finished by Theorem \ref{W1'}.
 \end{proof}

We illustrate Corollary \ref{C1'} by a typical example  of ergodic diffusion process which is not exponential ergodic.

 \beg{exa} Consider the following diffusion process $X_t$ on $\R^d$ generated by the operator $L$:
 $$Lf(x):=\DD f(x) - \<\nn V(x),\nn f(x)\>,\ \ \ x\in\R^d,$$
where  $V\in C^2(\R^d)$ satisfying
$$\sup_{x\in\R^d}|V(x)+\aa |x|^\bb|<\infty$$
for some constants $\aa>0$ and $\bb\in (0,1)$. Then $M_q(\mu)<\infty$ for any $q\in [1,\infty)$.

According to \cite[Example 1.4]{RW01} and the stability of the  weak  Poincar\'e inequality by bounded perturbations of $V$,   there exist constants $c_1,c_2\in (0,\infty)$ such that \eqref{XI} holds for
$$\xi(t)= c_1 \exp\Big[-c_2 t^{\ff{\bb}{4-3\bb}}\Big],\ \ \ \ t>0.$$
So, by Corollary \ref{C1'}, for any $p\in [1,\infty)$ there exists a constant
$c\in (0,\infty)$ such that
\beg{align*}&\sup_{\ll\in \scr P_0(\S)} \E^\mu\big[\W_{p,\ll}(\mu_t,\mu)^{2p}\big]\le \ff c t,\ \ \ t>0,\\
&\sup_{\ll\in \scr P_0(\S)} \E^\mu\big[\W_{p,\ll}(\tt\mu_n,\mu)^{2p}\big]\le \ff c n,\ \ \ n\in\mathbb N,\end{align*}
where, according to \eqref{ADS}, the orders $t^{-1}$ and $n^{-1}$  are sharp when $p=1$.
  \end{exa}

\section{Application to partially dissipative SPDEs  }

Since our results provide dimension-free convergence rates of empirical measures for ergodic Markov processes, they apply to infinite-dimensional models.
In this section,  we consider semi-linear SPDEs with long distance dissipation.

Let $(\H,\<\cdot,\cdot\>,\|\cdot\|)$   be a separable Hilbert space,   let    $(A,\D(A))$ be  a positive definite self-adjoint operator such that   $\{S_t=\e^{-At}\}_{ t\ge 0}$    is  a
$C_0$-contraction semigroup
on $\H$, and let
$\L_b(\H)$ be the space of bounded linear operators on $\H$ equipped with the operator norm $\|\cdot\|_{op}.$ Let
$$b:   \H\to  \H,\ \ \si:    \H\to \L_b(\H)$$ be   measurable maps.
 Consider the following SPDE on $\H$:
\beq\label{3*E} \d X_t= \big\{  b(X_t)- AX_t \big\}\d t +\si(X_t)\d W_t,\ \ t\in [0,T],\end{equation}
where $W_t$ is a cylindrical Brownian motion on $\H$ with respect to a complete filtered probability space  $(\OO, \{\F_t\}_{t\ge 0},\P)$, i.e. for any orthonormal family $\{h_i\}_{i=1}^n\subset \H$,
$\{\<h_i,W_t\>\}_{i=1}^n$ are independent one-dimensional Brownian motions.

  An $\H$-valued progressively measurable process $(X_t)_{t\ge 0}$ is called a mild solution to $(\ref{3*E})$, if   $\P$-a.s. for any $t\ge 0$,
\beg{align*}&  \int_0^t\big(|S_{t-s}b( X_s)|+\|S_{t-s} \si( X_{s})\|_{HS}^2\big)\d s<\infty,\\
&X_t = S_t X_0 + \int_0^t S_{t-s} b(s,X_s)\d s +\int_0^tS_{t-s}\si(s,X_s)\d W_s,\end{align*}
where $\|\cdot\|_{HS}$ is the Hilbert-Schmidt norm.

\beg{enumerate} \item[$(B_1)$] The maps $\si$  and $b$ are Lipschitz continuous, $\si$ is bounded, and  $ A$ has discrete spectrum with eigenvalues $\{\aa_i>0\}_{i\ge 1}$ listed in the increasing order counting multiplicities satisfying
 $$\sum_{i=1}^\infty \ff 1{\aa_i}<\infty.$$
 \item[$(B_2)$]  $\si\si^*= \kk^2 I + \hat\si\hat\si^*$ for some  $\kk\in (0,\infty)$ and Lipschitz continuous $\hat \si: \H\to \L_b(\H)$, where $I$ is the identity operator on $\H$. Moreover,
 there exist constants $ K_1, K_2, l\in (0,\infty) $    such that
\beq\label{DSL}\beg{split}&\ff 1 2 \|\hat\si(h_1)-\hat\si(h_2)\|_{HS}^2 + \<h_1-h_2, b(h_1)-b(h_2)- A(h_1-h_2)\>\\
& \le \big(K_1 1_{\{\|h_1-h_2\|\le l\}}- K_2 1_{\{\|h_1-h_2\|>l\}}\big)  \|h_1-h_2\|^2,\ \ h_1,h_2\in \H.\end{split}  \end{equation}
 \end{enumerate}

 It is well-known that condition $(B_1)$ implies  the existence of uniqueness of the mild solution $X_t^x$ for any initial value $x\in \H$, see \cite{DZ} or \cite[Theorem 3.1.1]{Wbook} for a more general result.
Consider  the associated Markov semigroup
 $$P_t f(x)= \E[f(X_t^x)],\ \ \ t\ge 0,\ x\in\H,\ f\in\B_b(\H).$$
 In general, let $P_t^*\nu$ be the distribution of the solution $X_t$ with initial distribution $\nu$. We have the following result.

  \beg{thm} Assume $(B_1)$ and $(B_2)$. Then the following assertions hold.
  \beg{enumerate} \item[$(1)$] There exist constants $c_0,\kk_0\in (0,\infty)$ such that
  $$\W_1(P_t^*\nu_1,P_t^*\nu_2)\le c_0\e^{-\kk_0 t} \W_1(\nu_1,\nu_2),\ \ \ t\ge 0,\ \nu_1,\nu_2\in \scr P(\H).$$
  \item[$(2)$] $P_t$ has a unique invariant probability measure $\mu$, which satisfies    $M_q(\mu)<\infty$ for any $q\in [1,\infty)$.
  \item[$(3)$]  There  exists a constant $c_0'\in (0,\infty)$ such that   $\eqref{XI}$ holds for
  $$\xi(t)= c_0' \e^{-\kk_0 t},\ \ \ t\ge 0.$$
   Consequently, for any $p\in [1,\infty)$ there exists a constant $c\in (0,\infty)$ such that for any $\ll\in \scr P_0(\S)$, $t>0$ and $n\in \mathbb N$,
   \beg{align*}&  \E^\mu\big[\W_{p,\ll}(\mu_t,\mu)^{2p}\big]\le   \ff c t,\ \   \
   \E^\mu\big[\W_{p,\ll}(\tt\mu_n,\mu)^{2p}\big]\le \ff c   n.\end{align*}
  \item[$(4)$]   There exists a constant $c'\in (0,\infty)$ such that for any $\nu\in \scr P_1(\H)$ and $\ll\in \scr P_0(\S)$,
    \beg{align*}&  \E^\nu\big[\W_{1,\ll}(\mu_t,\mu) \big]\le  c'\big(1+  M_1(\nu)\big) t^{-\ff 1 2},\ \  \ \ t>0,\\
 &  \E^\nu\big[\W_{1,\ll}(\tt\mu_n,\mu) \big]\le c'\big(1+  M_1(\nu)\big) n^{-\ff 1 2},\ \ \ \ n\in\mathbb N.\end{align*}
 \end{enumerate}\end{thm}

  \beg{proof}     Under $(B_2)$, we may reformulate \eqref{3*E} as
\beq\label{3*E'} \d X_t= \big\{  b(X_t)- AX_t \big\}\d t +\kk \d W_t^1 + \hat \si(X_t)\d W_t^2,\ \ t\in [0,T],\end{equation}
where $W_t^1$ and $W_t^2$ are independent cylindrical Brownian motions on $\H$.

(1)  Let $\{e_n\}_{n\ge 1}$ be the unitary eigenvectors of $A$ with respect to $\{\aa_n\}_{n\ge 1}$.
For any $n\in \mathbb N$, let $\H_n:={\rm span}\{e_i: \le i\le n\}$ and
$$\Pi_n h:= \sum_{i=1}^n \<h,e_i\> e_i,\ \ \ h\in\H$$ be the projection operator form $\H$ to $\H_n.$
For any $x \in\H $, let $ X_t^n $  solve the following SDE  on $\H_n$ with initial value $ x_n:= \Pi_n x$:
\beq\label{CP1} \d X_t^{x,n}= \big\{\Pi_n b(X_t^{x,n})-AX_t^{x,n}\big\} + \kk \Pi_n \d W_t^1+ \Pi_n \hat \si(X_t^{x,n}) \d W_t^2,\ \ \ X_0^{x,n}=x_n.\end{equation}
For any $x\ne y\in \H$, consider the following SDE for $X_0^{y,n}=y_n:=\Pi_n y$:
\beq\label{CP3}\beg{split}  &\d X_t^{y,n}= \big\{\Pi_n b(x_t^{y,n})-AX_t^{y,n}\big\} \\
&+ \kk \bigg(I- 21_{\{\tau>t\}}\ff{(X_t^{x,n}-X_t^{y,n})\otimes (X_t^{x,n}-X_t^{y,n})}{\|X_t^{x,n}-X_t^{y,n}\|^2}\bigg)\Pi_n\d W_t^1+ \Pi_n\hat \si(X_t^{y,n}) \d W_t^2,\end{split}\end{equation}
where $\tau:=\inf\{t\ge 0: X_t^n=Y_t^n\}$ is the coupling time.
By $(B_1)$ and $(B_2)$, this SDE have a unique solution with $X_t^{y,n}=X_t^{x,n}$ for $t\ge \tau$.

Let $ \gg(r):= K_1\big(r\land\ff{l^2}{r}\big)-K_2r,\  r>0.$ Then  \eqref{DSL} implies
\beq\label{DSL'}\beg{split}&\ff 1 2 \|\hat\si(h_1)-\hat\si(h_2)\|_{HS}^2 + \<h_1-h_2, b(h_1)-b(h_2)- A(h_1-h_2)\>\\
& \le \gg(\|h_1-h_2\|)  \|h_1-h_2\|,\ \ h_1,h_2\in \H.\end{split}  \end{equation}
 According to step (a) in the proof of \cite[Corollary 3.2]{W23},
the function
$$g(r)= \int_0^r \e^{-\ff {\gg(s)} {2\kk^2} }\d s\int_s^\infty t \e^{\ff{\gg(t)}{2\kk^2}}\d t,\ \ r\ge 0 $$
satisfies   \beq\label{QA'} g''(r)\le 0,\ \ \ \ \ r>0,\end{equation}
so that  $g(\|h_1-h_2\|)$ is a distance on $\H$.
Moreover, it is easy to see that
\beq\label{QA0'}  2\kk^2 g''(r)+ \gg(r)g'(r)= - r,\ \ \ r>0,\end{equation}
\beq\label{QA1'} \bb_1:= \inf_{r>0} \ff r{g(r)}>0,\ \  \ \bb_2:=\sup_{r>0} \ff r{g(r)}<\infty. \end{equation}
 By \eqref{CP1}, \eqref{CP3}, \eqref{DSL'} and It\^o's formula, we obtain
\beg{align*} &\d \|X_t^{x,n}-X_t^{y,n}\|-  \d M_t\\
&\le \bigg\{\Big\<\ff{X_t^{x,n}-X_t^{y,n}}{\|X_t^{x,n}-X_t^{y,n}\|},\ b(X_t^{x,n})-b(X_t^{y,n})- A(X_t^{x,n}-X_t^{y,n})\Big\>_\H+ \ff {\| \hat \si(X_t^{x,n})-\hat\si(X_t^{y,n}) \|_{HS}^2}{2\|X_t^{x,n}-X_t^{y,n}\|}
\bigg\}\d t\\
&\le   \gg(\|X_t^{x,n}-X_t^{y,n}\|)\d t,\ \ \ t<\tau,\end{align*}
where
$$\d M_t= \bigg\<\ff{X_t^{x,n}-X_t^{y,n}}{\|X_t^{x,n}-X_t^{y,n}\|}, \ 2\kk \d W_t^1+ \Pi_n \{\hat\si(X_t^{x,n})-\hat\si(X_t^{y,n})\}\d W_t^2\bigg\>_\H $$ is a martingale with $\d\<M\>_t\ge 4\kk^2\d t.$
This together with  \eqref{QA'},  \eqref{QA0'}  and \eqref{QA1'}  yields
\beg{align*}&\d g(\|X_t^{x,n}-X_t^{y,n}\|)- g'(\|X_t^{x,n}-X_t^{y,n}\|)\d M_t \\
&\le \Big\{\gg(\|X_t^{x,n}-X_t^{y,n}\|) g'(\|X_t^{x,n}-X_t^{y,n}\|)+2\kk^2 g''(\|X_t^{x,n}-X_t^{y,n}\|)\Big\}\d t\\
&= -\|X_t^{x,n}-X_t^{y,n}\|\d t\le - \bb_1 g(\|X_t^{x,n}-X_t^{y,n}\|)\d t,\ \ \ t<\tau.\end{align*}
Noting that   $X_t^{x,n}=X_t^{y,n}$ for $t\ge \tau$, this together with \eqref{QA1'} implies
$$\E [\|X_t^{x,n}-X_t^{y,n}\| ] \le \bb_2 \E \big[g(\|X_t^{x,n}-X_t^{y,n}\|)\big]\le \bb_2 \e^{-\bb_1 t} g(\|x_n-y_n\|)\le \ff{\bb_2}{\bb_1} \e^{-\bb_1 t} \|x_n-y_n\|,\ \ \ t\ge 0.$$
Combining this with \cite[Theorem 3.1.2]{Wbook}, by letting $n\to\infty$ we find a coupling $(X_t^x,X_t^y)$ for the Markov process with semigroup $P_t$ stating from $(x,y)$ such that
\beq\label{RRW}   \W_1(P_t^*\dd_x,P_t^*\dd_y) \le \E [\|X_t^{x}-X_t^{y}\| ]\le  \ff{\bb_2}{\bb_1} \e^{-\bb_1 t} \|x-y\|,\ \ \ t\ge 0.\end{equation}
This implies the desired estimate on $ \W_1(P_t^*\nu_1,P_t^*\nu_2)$ for $c_0:= \ff{\bb_2}{\bb_1}$ and $\kk_0=\bb_1.$

(2)  By the estimate in (1)  and \cite[Theorem 1.6.4]{WR25}, for the existence and uniqueness of invariant probability measure, we only need to show that
$M_1(P_t^*\nu)$  is locally bounded in $t$ for any fixed $\nu\in \scr P_1(\H)$.

Let $X_0$ satisy $\nu:=\L_{X_0}\in \scr P_1(\H)$, and let
$$T_n:= \inf\big\{t\ge 0: \|X_t\|\ge n\big\},\ \ \ n\ge 1.$$
We have $T_n\to \infty$ as $n\to\infty$.
By $(A_1)$, we find   constants $c_1,c_2\in (0,\infty)$ such that
\beg{align*} & \E[\|X_{t\land T_n}\|]\le \E[\|X_0\|] + c_1 \int_0^t \big(1+  \E[\|X_{s\land T_n}\|]\big)\d s + c_1 \bigg(\int_0^t \|S(t-s)\|_{HS}^2\d s\bigg)^{\ff 1 2}\\
&\le M_1(\nu) + c_2 t +c_1 \int_0^t  \E[\|X_{s\land T_n}\|] \d s,\ \ \ t\ge 0.\end{align*}
By Gronwall's inequality, we obtain
$$\E[\|X_{t\land T_n}\|] \le \big[M_1(\nu) + c_2 t \big] \e^{c_1 t},\ \ \ t\ge 0, \ n\ge 1.$$
By Fatou's lemma with $n\to\infty$, we derive the local boundedness of  $M_1(P_t^*\nu)= \E[\|X_t\|]$.

To prove $M_q(\mu)<\infty$ for any $q\ge 2$, let $X_0= 0$ and  $Y_t$ solve the equation \eqref{3*E} with $b=0$, i.e.
$$Y_t=  \int_0^t S(t-s)   \si(Y_s)\d W_s^2,\ \ \ t\ge 0.$$
By $(A_1)$ we find a constant $c_3\in (0,\infty)$ such that
\beq\label{QY} \E [\|Y_t\|^q ]\le   c_3 \|_\infty \bigg(\int_0^t \|S(t-s)\|_{HS}^2\d s\bigg)^{\ff q 2}\le c_3\bigg(\sum_{i=1}^\infty \ff 1 {\aa_i}\bigg)^{\ff q 2}<\infty.\end{equation}
Next, by $(A_1)$ and It\^o's formula for $X_t-Y_t$, we find a constant $c_4 \in (0,\infty)$ such that
\beg{align*} &\d \|X_t-Y_t\|^2\le 2\big(K_1 l^2 - K_2 \|X_t-Y_t\|^2+ 2 \|b(Y_t)\| \|X_t-Y_t\| \big)\d t +\d M_t\\
&\le \big[c_4 (1+\|Y_t\|^4) -  K_2 \|X_t-Y_t\|^2\big]\d t +\d M_t,\end{align*}
where $M_t:= 2\<X_t-Y_t, (\hat \si(X_t)-\hat \si(Y_t))\d W_t^2\>$ satisfies
$$\d\<M\>_t\le c_4 \|X_t-Y_t\|^2\d t.$$
So, for any $q\ge 2$, there exists a constant $c(q)\in (q,\infty)$ such that
$$\d \|X_t-Y_t\|^q\le \Big[c(q) (1+\|Y_t\|^{c(q)}) - \ff {K_2}2  \|X_t-Y_t\|^q\Big]\d t + \ff q 2 \|X_t-Y_t\|^{q-2} \d M_t,$$
which, together with Gronwall's inequality, implies
$$\sup_{t\ge 0}\E[\|X_t-Y_t\|^q]<\infty.$$
This and \eqref{QY} implies $\sup_{t\ge 0} \E[\|X_t\|^q]<\infty$.     Since $X_t$ converges to $\mu$ weakly as $t\to\infty$, we derive $M_q(\mu)<\infty$.

(3) By Corollary \ref{C1'} and (2), it suffices to verify \eqref{XI} for the claimed $\xi(t)$.

For any bounded  Lipschitz continuous function $f$ on $\H$, let $\|\nn f\|_\infty$ be the Lipschitz constant. By \eqref{RRW}, we have
\beq\label{CP4}\beg{split}&|P_t f(x )-P_t f(y )| =    |\E[f(X_t^{x })- f(X_t^{y })]|\\
&\le \|\nn f\|_\infty  \E[\|X_t^{x }-X_t^{y }\| ] \le \ff{\bb_2\|\nn f\|_\infty}{\bb_1} \e^{-\bb_1 t} \|x -y \|.\end{split} \end{equation}
Next, by \cite[Theorem 3.3.1]{Wbook}, there exists a constant $k_1\in (0,\infty)$ such that
$$\|\nn P_1 f\|_\infty \le k_1\|f\|_\infty,\ \ \ f\in \B_b(\H).$$ Combining this with  \eqref{CP4},
if the invariant probability measure $\mu$ satisfies $M_2(\mu)<\infty$, we find a constant $k_2\in (1,\infty)$ such that
 $$\|P_{t+1}f -\mu(f)\|_{L^2(\mu)}^2 = \ff 1 2 \int_{\H\times \H} |P_{t+1}f(x)-P_{t+1}f(y)|^2 \mu(\d x)\mu(\d y) \le k_2^2 \e^{-2\bb_1 t} \|f\|_\infty,\ \ f\in \B_b(\H).$$
Since $\|P_t-\mu\|_{L^\infty(\mu)\to L^2(\mu)} \le 1$, we conclude that \eqref{XI} holds for
$$\xi(t)= k_2\e^{\bb_1} \e^{-\bb_1 t},\ \ \ t\ge 0.$$

(4) By \eqref{RRW} and Proposition \ref{PR} for $p=1$,   the desired estimate follows from that in (3).
\end{proof}

\paragraph{Acknowledgements.} The author would like to thank Xiaochen Ma and Panpan Ren for corrections.

\paragraph{Data availability statement.}
No new data were created or analysed in this study. Data sharing is not applicable.
\paragraph{Conflict of interest statement.}
The author  has no conflicts of interest to declare.

\end{document}